\documentclass[11pt,reqno]{amsart}

\usepackage{amsmath,amssymb,amsthm}
\usepackage[margin=1in]{geometry}
\usepackage{mathtools}
\usepackage[dvipsnames]{xcolor}
\usepackage[colorlinks=true,linkcolor=lava,citecolor=blue,urlcolor=brown]{hyperref}
\definecolor{lava}{rgb}{0.8,0.2,0.0}
\usepackage{orcidlink}

\theoremstyle{plain}
\newtheorem{theorem}{Theorem}[section]

\newtheorem{lemma}{Lemma}[section]
\newtheorem{corollary}{Corollary}[section]

\theoremstyle{definition}
\newtheorem{definition}{Definition}[section]
\newtheorem{example}{Example}[section]
\newtheorem{problem}{Open problem}[section]

\theoremstyle{remark}
\newtheorem{remark}{Remark}[section]

\numberwithin{equation}{section}

\newcommand{\RN}{\mathbb{R}^N}
\newcommand{\R}{\mathbb{R}}
\newcommand{\dx}{\,dx}
\newcommand{\dt}{\,dt}
\newcommand{\ds}{\,ds}
\newcommand{\eps}{\varepsilon}

\newcommand{\Q}{Q}
\newcommand{\pF}{p_{F}}
\DeclareMathOperator{\supp}{supp}
\newcommand{\abs}[1]{\left| #1 \right|}
\newcommand{\norm}[1]{\left\| #1 \right\|}
\newcommand{\iQ}{\iint_{\Q}}
\newcommand{\tn}[1]{[\![ #1 ]\!]}          
\newcommand{\Duh}{\mathcal{D}}             
\newcommand{\Lop}{\mathcal{L}}             
\newcommand{\hN}{h_{N}}
\newcommand{\conv}{\star}
\newcommand{\Meff}{M}

\begin{document}

\title[Time-integrable forcing and the Fujita phenomenon]{The Fujita equation with time-integrable forcing:\\ effective mass, mass annihilation and sharp lifespan}

\author[ Ahmad Z. Fino and M. Majdoub]{Ahmad Z. Fino and Mohamed Majdoub}

\address[Ahmad Z. Fino]{Department of Mathematics, Faculty of Sciences, Lebanese University, P.O. Box 1352 Tripoli, Lebanon.}
\email{afino@ul.edu.lb}
\address[M. Majdoub]{Department of Mathematics, College of Science, Imam Abdulrahman Bin Faisal University, P. O. Box 1982, Dammam, Saudi Arabia.}
\address[M. Majdoub]{Basic and Applied Scientific Research Center, Imam Abdulrahman Bin Faisal University, P.O. Box 1982, 31441, Dammam, Saudi Arabia.}
\email{mmajdoub@iau.edu.sa}
\email{med.majdoub@gmail.com}
\email{mohamed.majdoub@fst.rnu.tn}


\date{\today}

\begin{abstract}
We study the Cauchy problem $\partial_t u-\Delta u=\abs{u}^{p}+\zeta(t)\,w(x)$ in $\RN\times(0,\infty)$, $p>1$, with sign-changing initial data and a forcing term whose temporal profile $\zeta$ is integrable on $(0,\infty)$; the model case is $\zeta(t)=t^{\sigma}e^{-t}$, $\sigma>-1$. For power-type profiles $\zeta(t)=t^{\sigma}$ the forcing is known either to raise the critical exponent above the Fujita exponent $\pF=1+2/N$ or to destroy global existence altogether. We show that an integrable forcing leaves $\pF$ unchanged and acts only through the \emph{effective mass}
$\Meff=\int_{\RN}u_0\dx+\big(\int_0^\infty\zeta\dt\big)\big(\int_{\RN}w\dx\big)$.
Our central result is an exact \emph{mass-annihilation identity}: if $1<p\le\pF$, every global weak solution lies in $L^{p}(\RN\times(0,\infty))$ and satisfies $\iint\abs{u}^{p}\dx\dt=-\Meff$. Consequently, in the Fujita range there is no global solution as soon as $\Meff\ge0$ and the data are nontrivial (the borderline case $\Meff=0$ included), and every global mild solution vanishes in $L^1(\RN)$ as $t\to\infty$; nonpositive data show that $\Meff<0$ is compatible with global existence, while large localized bumps show that it does not imply it. For $p>\pF$ we prove that every solution is dominated by a fixed multiple of the linear forced flow as soon as a single integral condition holds; this gives small-data global existence without sign restrictions and, under a mild decay assumption on $\zeta$, convergence to the profile $m_\infty G_t$ with the explicitly identified asymptotic mass $m_\infty=\Meff+\iint\abs{u}^{p}$. We complement these results with a localized necessary condition valid for every $p>1$, which shows that large or slowly decaying data blow up regardless of the forcing and that non-integrable forcings shift the threshold, and with sharp lifespan estimates for small data with positive effective mass: $T_\eps\simeq\eps^{-2(p-1)/(2-N(p-1))}$ for $p<\pF$ and $\log T_\eps\simeq\eps^{-(p-1)}$ for $p=\pF$.
\end{abstract}

\maketitle


\section{Introduction}\label{sec:intro}

We consider the Cauchy problem for the semilinear heat equation with a source and a separable external forcing,
\begin{equation}\label{eq:model}
\begin{cases}
\partial_t u-\Delta u=\abs{u}^{p}+\zeta(t)\,w(x), & (x,t)\in \RN\times(0,\infty),\\[2pt]
u(\cdot,0)=u_0, & x\in\RN,
\end{cases}
\end{equation}
where $N\ge1$, $p>1$, and $u_0$, $w$ and $\zeta$ are real-valued and allowed to change sign. The distinguishing feature of the present work is that the temporal profile of the forcing is \emph{integrable},
\[
\zeta\in L^{1}(0,\infty),
\]
the model example being
\begin{equation}\label{eq:model-zeta}
\zeta(t)=t^{\sigma}e^{-t},\qquad \sigma>-1,\qquad \int_0^\infty\zeta\dt=\Gamma(\sigma+1).
\end{equation}
Throughout the paper, we write $\Q:=\RN\times(0,\infty)$, $p'=p/(p-1)$,
\[
\pF:=1+\frac{2}{N}
\]
for the Fujita exponent, and
\[
Z:=\int_0^{\infty}\zeta(t)\dt,\qquad W:=\int_{\RN}w(x)\dx,\qquad m_0:=\int_{\RN}u_0(x)\dx .
\]
The number
\begin{equation}\label{eq:total-impulse}
\Meff:=m_0+Z\,W
\;=\;\int_{\RN}u_0\dx+\Big(\int_0^{\infty}\zeta\dt\Big)\Big(\int_{\RN}w\dx\Big)
\end{equation}
will be called the \emph{effective mass} of the data: it is the bare mass of the initial datum augmented by the total impulse that the forcing injects over the whole time axis. In the model case \eqref{eq:model-zeta} one has $\Meff=m_0+\Gamma(\sigma+1)W$.

\subsection{Background: the Fujita phenomenon}
In his seminal paper \cite{Fujita1966}, Fujita showed that for $\partial_tu=\Delta u+u^{p}$ with nonnegative data every nontrivial solution blows up in finite time when $1<p<\pF$, whereas global positive solutions exist for small data when $p>\pF$. The critical case $p=\pF$ belongs to the blow-up range; this was proved by Hayakawa \cite{Hayakawa1973} for $N=1,2$, by Kobayashi, Sirao and Tanaka \cite{KobayashiSiraoTanaka1977} in all dimensions, and by different methods by Aronson and Weinberger \cite{AronsonWeinberger1978} and Weissler \cite{Weissler1981}. The local theory in Lebesgue spaces goes back to Weissler \cite{Weissler1980}; the eigenfunction method of Kaplan \cite{Kaplan1963} is one of the basic tools for blow-up. We refer to the surveys of Levine \cite{Levine1990} and Deng and Levine \cite{DengLevine2000} and to the monograph of Quittner and Souplet \cite{QuittnerSouplet2019} for the extensive literature that followed.

For the equation $\partial_tu-\Delta u=\abs{u}^{p}$ with sign-changing data, the nonnegativity of the source makes the problem amenable to the nonlinear capacity (test-function) method of Mitidieri and Pohozaev \cite{MitidieriPokhozaev2001}, which works directly with weak solutions and shows nonexistence of global solutions in the Fujita range under a sign condition on $\int u_0$. Sign-changing solutions exhibit many phenomena absent in the positive case; see for instance \cite{CazenaveDicksteinWeissler2009} and \cite[Chapter~II]{QuittnerSouplet2019}. The size of the lifespan of small solutions is a closely related question: for $u_0=\eps a$ with $a\ge0$, Lee and Ni \cite{LeeNi1992} proved that $T_\eps\simeq\eps^{-(1/(p-1)-N/2)^{-1}}$ for $p<\pF$ and $\log T_\eps\simeq\eps^{-(p-1)}$ for $p=\pF$, and Ikeda and Sobajima \cite{IkedaSobajima2019} devised a test-function method yielding sharp upper bounds for weak solutions of critical problems.

A second classical theme that will play a role below is the behavior of the mass. For the heat equation with absorption $\partial_tv=\Delta v-v^{p}$, $v_0\ge0$, the mass $\int v(t)$ converges to a positive limit when $p>\pF$ but tends to zero when $p\le\pF$; see Gmira and V\'eron \cite{GmiraVeron1984}, Kamin and Peletier \cite{KaminPeletier1985}, Brezis, Peletier and Terman \cite{BrezisPeletierTerman1986}, Galaktionov, Kurdyumov and Samarskii \cite{GalaktionovKurdyumovSamarskii1986}, Escobedo, Kavian and Matano \cite{EscobedoKavianMatano1995} and Herraiz \cite{Herraiz1999}. Analogous decay-of-mass results were obtained for the fractional Laplacian by Fino and Karch \cite{FinoKarch2010} and for mixed local--nonlocal operators by Kirane, Fino and Ayoub \cite{KiraneFinoAyoub2025}. One of the outcomes of the present paper is that the same threshold governs the mass of \emph{every} global solution of the source problem \eqref{eq:model}, forced or not.

\subsection{Forcing terms}
The influence of an external forcing on the Fujita phenomenon was first understood for stationary forcings. Bandle, Levine and Zhang \cite{BLZ2000} studied $\partial_tu=\Delta u+u^{p}+w(x)$ and proved that, when $\int w>0$, the critical exponent jumps from $\pF$ to $N/(N-2)$ if $N\ge3$, and that there is no global solution at all if $N\in\{1,2\}$; see also Zhang \cite{Zhang1999}. For space--time forcings $t^{\sigma}w(x)$, Jleli, Kawakami and Samet \cite{JKS2020} found the critical exponent
\[
p^{*}(\sigma)=\frac{N-2\sigma}{N-2\sigma-2}\quad(-1<\sigma<0),\qquad p^{*}(\sigma)=\infty\quad(\sigma>0),
\]
which is discontinuous at $\sigma=0$ when $N\ge3$. These results were extended to systems by Fino, Jleli and Samet \cite{FinoJleliSamet2020}, to inhomogeneous and Hardy--H\'enon nonlinearities in \cite{Majdoub2021,Majdoub2023}, to general time profiles $a(t)w(x)$ controlled by the averages $t^{-1}\int_0^ta$ in \cite{AlshehriEtAl2023}, to nonlocal nonlinearities and mixed local--nonlocal or space--time fractional diffusions by Ben Belgacem and Majdoub \cite{BenBelgacemMajdoub2025a,BenBelgacemMajdoub2025b,BenBelgacemMajdoub2025c}, to degenerate and singular weights in \cite{MajdoubTorebek2026}, and to spatial profiles with regularly varying mass in \cite{KumarMajdoub2026}. Related questions for time--space fractional problems were studied by Fino and Kirane \cite{FinoKirane2012}, and the robustness of $\pF$ under singular perturbations of the diffusion was established in \cite{MajdoubMliki2026,MajdoubTorebek2026b}.

In all these works the forcing is \emph{non-integrable in time} (the power $t^{\sigma}$ with $\sigma>-1$ has $\int_0^Tt^{\sigma}\dt\to\infty$), and its accumulated impulse $\int_0^T\zeta$ grows without bound. This is precisely what shifts the critical exponent: in the test-function balance the forcing contributes a term of size $T^{\sigma+1}$, which competes with the diffusive cost $R^{N+2-2p'}$ at the parabolic scale $T=R^{2}$. As $\sigma\downarrow-1$ one has $p^{*}(\sigma)\downarrow\pF$, which suggests that integrable forcings sit at an endpoint of the theory. It is natural to ask what exactly happens there.

\subsection{Questions and heuristics}
When $\zeta\in L^{1}(0,\infty)$ the forcing delivers a \emph{finite} total impulse $Z\,W$. Integrating \eqref{eq:model} in space, any solution with integrable data satisfies the mass balance
\begin{equation}\label{eq:mass-balance-intro}
\frac{d}{dt}\int_{\RN}u(t)\dx=\int_{\RN}\abs{u(t)}^{p}\dx+\zeta(t)\,W ,
\end{equation}
so that, formally, $\int u(t)$ tends to $\Meff+\iint_\Q\abs{u}^{p}$ as $t\to\infty$. This suggests that an integrable forcing acts as a mass shift and that the relevant quantity is the effective mass \eqref{eq:total-impulse}. The following questions arise naturally, and the purpose of this paper is to answer them.

\emph{Is $\pF$ still the critical exponent?} Since $\zeta$ is integrable, one may hope that small forcings do not destroy global existence above $\pF$, while the extra mass $ZW$ cannot prevent blow-up below $\pF$ when it is positive.

\emph{Which data blow up in the Fujita range?} For sign-changing data, the test-function method yields nonexistence when the relevant mass is \emph{positive}. What happens in the borderline case $\Meff=0$, and when $\Meff<0$?

\emph{What is the large-time behavior of global solutions, and what is the lifespan of small blowing-up solutions?}

\subsection{Main results}
We now state our main results; the precise notions of weak and mild solutions are given in Definitions~\ref{def:weak} and \ref{def:mild}. For $f\in L^1(\RN)\cap L^\infty(\RN)$ we write $\tn{f}:=\max\{\norm{f}_{L^1},\norm{f}_{L^\infty}\}$, $G_t(x)=(4\pi t)^{-N/2}e^{-\abs{x}^2/4t}$ denotes the heat kernel, $S(t)f=G_t*f$, and
\begin{equation}\label{eq:hN}
\hN(t):=\min\{1,(4\pi t)^{-N/2}\},\qquad (\abs{\zeta}\conv \hN)(t):=\int_0^{t}\abs{\zeta(s)}\,\hN(t-s)\ds ,
\end{equation}
so that $\norm{S(t)f}_{L^\infty}\le \hN(t)\tn{f}$ for all $t>0$.

\subsubsection*{The Fujita range: mass annihilation}
Our central result is an exact identity satisfied by every global solution in the Fujita range.

\begin{theorem}\label{thm:annihilation}
Let $1<p\le\pF$, $\zeta\in L^1(0,\infty)$, $w\in L^1(\RN)$ and $u_0\in L^1(\RN)$. If $u$ is a global weak solution of \eqref{eq:model}, then $u\in L^{p}(\Q)$ and
\begin{equation}\label{eq:annihilation}
\iQ\abs{u}^{p}\dx\dt=-\Meff=-\Big(\int_{\RN}u_0\dx+Z\,W\Big).
\end{equation}
\end{theorem}

The identity \eqref{eq:annihilation} says that the total amount of mass produced by the source over the whole time axis exactly cancels the effective mass of the data. Its first consequence is a sharp nonexistence statement, which replaces the sign condition $\Meff>0$ by $\Meff\ge0$.

\begin{theorem}\label{thm:nonexist}
Assume $1<p\le\pF$, $\zeta\in L^1(0,\infty)$, $w\in L^1(\RN)$ and $u_0\in L^1(\RN)$.
\begin{enumerate}
\item[(i)] If
\begin{equation}\label{eq:mass-condition}
\Meff=m_0+Z\,W>0 ,
\end{equation}
then \eqref{eq:model} admits no global weak solution.
\item[(ii)] If $\Meff=0$, then \eqref{eq:model} admits a global weak solution if and only if $u_0=0$ a.e. and $\zeta(t)w(x)=0$ for a.e. $(x,t)\in\Q$; in that case $u\equiv0$ is the only global weak solution.
\item[(iii)] If $\Meff<0$, every global weak solution satisfies $\norm{u}_{L^p(\Q)}^{p}=\abs{\Meff}$.
\end{enumerate}
\end{theorem}

For the unforced equation $\partial_tu-\Delta u=\abs{u}^p$, statements (i)--(ii) are in line with the nonexistence results obtained for the corresponding differential inequality by the nonlinear capacity method \cite{MitidieriPokhozaev2001}: nontrivial data with zero mass, such as $u_0=a-a(\cdot-x_1)$, never give rise to global solutions when $p\le\pF$. What is new is the exact identity \eqref{eq:annihilation}, its consequences for the dynamics of global solutions, and the fact that in the forced problem the relevant quantity is the effective mass $\Meff$ and not $m_0$. Part (ii) shows, for instance, that the datum $u_0=-Z\,w$ (zero effective mass) blows up whenever $\zeta w\not\equiv0$. The next result shows that both alternatives genuinely occur when $\Meff<0$.

\begin{theorem}\label{thm:nonpositive}
Let $p>1$, $\zeta\in L^1(0,\infty)$, and let $u_0,w\in L^1(\RN)\cap L^\infty(\RN)$ satisfy $u_0\le0$ and $\zeta(t)w(x)\le0$ for a.e. $(x,t)\in\Q$. Then \eqref{eq:model} has a unique global mild solution, and
\[
-\Big(S(t)\abs{u_0}+\int_0^tS(t-s)\abs{\zeta(s)w}\ds\Big)\le u(t)\le0\qquad\text{for all }t>0 .
\]
If moreover $1<p\le\pF$, then $\iQ\abs{u}^{p}=\norm{u_0}_{L^1}+\norm{\zeta w}_{L^1(\Q)}$ and $\norm{u(t)}_{L^1}\to0$ as $t\to\infty$.
\end{theorem}

On the other hand, data with $\Meff<0$ containing a large positive bump blow up for every $p>1$ (Corollary~\ref{cor:obstruction}(a) and Example~\ref{ex:both}). Thus, in the Fujita range, the sign of the effective mass decides the alternative completely when $\Meff\ge0$, while for $\Meff<0$ both alternatives occur. This answers the question of whether an arbitrarily negative effective mass is still sufficient to force blow-up.

\subsubsection*{Large-time behavior and the asymptotic mass}
The identity \eqref{eq:annihilation} has a striking dynamical counterpart.

\begin{theorem}\label{thm:asymptotics}
Let $p>1$, $\zeta\in L^1(0,\infty)$ and $u_0,w\in L^1(\RN)\cap L^\infty(\RN)$.
\begin{enumerate}
\item[(i)] If $u$ is a global mild solution with $\abs{u}^p\in L^1(\Q)$, then, with $m_\infty:=\Meff+\iQ\abs{u}^p$,
\[
\lim_{t\to\infty}\norm{u(t)-m_\infty G_t}_{L^1}=0 .
\]
\item[(ii)] If $1<p\le\pF$, every global mild solution satisfies the hypothesis of (i) with $m_\infty=0$; in particular $\norm{u(t)}_{L^1}\to0$ as $t\to\infty$.
\item[(iii)] If $p>\pF$ and $\zeta$ satisfies hypothesis \textup{(H$_a$)} below, then every global mild solution satisfying $\abs{u}\le\kappa V$ (see Theorem~\ref{thm:exist}) has $\abs{u}^p\in L^1(\Q)$ and, for every $q\in[1,\infty]$,
\[
\lim_{t\to\infty}t^{\frac N2(1-\frac1q)}\norm{u(t)-m_\infty G_t}_{L^q}=0 ,
\qquad m_\infty\ge \Meff .
\]
\end{enumerate}
\end{theorem}

Here and below we use the decay hypothesis
\begin{itemize}
\item[(H$_a$)] $\zeta\in L^1(0,\infty)$ and there exist $C_\zeta>0$ and $a>\max\{1,N/2\}$ such that $\abs{\zeta(t)}\le C_\zeta\,t^{-a}$ for a.e.\ $t\ge1$,
\end{itemize}
which holds for the model \eqref{eq:model-zeta} with every $a>0$. Theorem~\ref{thm:asymptotics} yields a characterization of the Fujita exponent through the asymptotic mass (Corollary~\ref{cor:characterization}): $p\le\pF$ if and only if every global solution loses all its mass.

\subsubsection*{Global existence above the Fujita exponent}
Global existence is obtained through a comparison with the \emph{linear forced flow with absolute data}
\begin{equation}\label{eq:V-intro}
\begin{aligned}
V(x,t)&:=\big[S(t)\abs{u_0}\big](x)+\int_0^{t}\abs{\zeta(s)}\,\big[S(t-s)\abs{w}\big](x)\ds,\\
m(t)&:=\norm{S(t)\abs{u_0}}_{L^\infty}+\int_0^{t}\abs{\zeta(s)}\,\norm{S(t-s)\abs{w}}_{L^\infty}\ds .
\end{aligned}
\end{equation}

\begin{theorem}\label{thm:exist}
Let $p>1$, $u_0,w\in L^1(\RN)\cap L^\infty(\RN)$, $\zeta\in L^1_{\rm loc}([0,\infty))$ and $0<T\le\infty$. Set $\kappa:=p/(p-1)$ and $c_p:=(p-1)^{p-1}/p^{p}$. If
\begin{equation}\label{eq:criterion}
\int_0^{T}m(t)^{p-1}\dt\le c_p ,
\end{equation}
then \eqref{eq:model} has a unique mild solution on $(0,T)$, and it satisfies the pointwise bound
\begin{equation}\label{eq:domination}
\abs{u(x,t)}\le\kappa\,V(x,t)\qquad\text{for a.e. }(x,t)\in\RN\times(0,T).
\end{equation}
No sign condition on $u_0$, $w$ or $\zeta$ is required.
\end{theorem}

When $T=\infty$, \eqref{eq:criterion} can only hold if $p>\pF$ (unless the data vanish), and it is then a smallness condition. Writing $\omega_\zeta:=\int_0^\infty(\abs{\zeta}\conv\hN)(t)^{p-1}\dt$, we obtain the following.

\begin{corollary}\label{cor:small}
Let $p>\pF$ and $\zeta\in L^1_{\rm loc}([0,\infty))$ with $\omega_\zeta<\infty$. There exists $\delta=\delta(N,p)>0$ such that, if $u_0,w\in L^1(\RN)\cap L^\infty(\RN)$ satisfy
\[
\tn{u_0}^{p-1}+\omega_\zeta\,\tn{w}^{p-1}\le\delta ,
\]
then \eqref{eq:model} has a unique global mild solution, $\abs{u}\le\kappa V$, and for all $t>0$
\[
\norm{u(t)}_{L^\infty}\le\kappa\big(\tn{u_0}\hN(t)+\tn{w}(\abs{\zeta}\conv\hN)(t)\big),\qquad
\norm{u(t)}_{L^1}\le\kappa\Big(\norm{u_0}_{L^1}+\norm{w}_{L^1}\int_0^t\abs{\zeta}\Big).
\]
Moreover $\omega_\zeta<\infty$ in each of the following cases: \textup{(a)} $p\ge2$ and $\zeta\in L^1(0,\infty)$; \textup{(b)} $\zeta$ satisfies \textup{(H$_a$)}, in which case $\norm{u(t)}_{L^\infty}\le C(1+t)^{-N/2}$; \textup{(c)} $N\ge3$, $\zeta\in L^1(0,\infty)$ and $\abs{\zeta(t)}\le Ct^{-b}$ for $t\ge1$ with $b>1$ and $b(p-1)>1$.
\end{corollary}

In particular, for the model \eqref{eq:model-zeta}, small data of either sign give global solutions for every $p>\pF$, with $\norm{u(t)}_{L^\infty}\le C\eps(1+t)^{-N/2}$. Hence the critical exponent for integrable forcings is exactly $\pF$; nonexistence for every $p>1$, which holds for growing forcings $t^{\sigma}w$, $\sigma>0$ \cite{JKS2020}, fails here.

\subsubsection*{A localized obstruction valid for all $p>1$}
The test-function method also produces a necessary condition for global existence that does not use the integrability of the data.

\begin{theorem}\label{thm:obstruction}
Let $p>1$, $u_0\in L^1_{\rm loc}(\RN)$, $w\in L^1_{\rm loc}(\RN)$ and $\zeta\in L^1_{\rm loc}([0,\infty))$, and let $\eta$ be the cutoff fixed in \eqref{eq:eta}. There is $C_0=C_0(N,p,\eta)>0$ such that, if \eqref{eq:model} has a global weak solution, then for all $x_0\in\RN$ and $R,T>0$,
\begin{multline}\label{eq:obstruction}
\int_{\RN}u_0\,\eta\Big(\tfrac{\abs{x-x_0}^2}{R^2}\Big)^{2p'}\dx
\\+\Big(\int_0^\infty\zeta(t)\,\eta\big(\tfrac tT\big)^{2p'}\dt\Big)\Big(\int_{\RN}w\,\eta\Big(\tfrac{\abs{x-x_0}^2}{R^2}\Big)^{2p'}\dx\Big)
\le C_0\big(R^{N}T^{1-p'}+R^{N-2p'}T\big).
\end{multline}
\end{theorem}

Corollary~\ref{cor:obstruction} derives from \eqref{eq:obstruction} that large data blow up for every $p>1$ whatever the forcing, that large forcings blow up for every $p>1$, and that data decaying like $\abs{x}^{-\gamma}$ with $\gamma<\min\{N,2/(p-1)\}$ blow up (the forced analogue of a result of Lee and Ni \cite{LeeNi1992}). Corollary~\ref{cor:nonintegrable} treats non-integrable forcings and recovers the blow-up part of \cite{JKS2020}.

\subsubsection*{Sharp lifespan}
Finally, we determine the lifespan of small solutions with positive effective mass in the Fujita range. For $a,b\in L^1(\RN)\cap L^\infty(\RN)$ and $\eps>0$, let $T_\eps\in(0,\infty]$ be the supremum of the $T>0$ for which \eqref{eq:model} with $u_0=\eps a$ and $w=\eps b$ has a mild solution on $(0,T)$.

\begin{theorem}\label{thm:lifespan}
Let $1<p\le\pF$, $\zeta\in L^1(0,\infty)$, $a,b\in L^1(\RN)\cap L^\infty(\RN)$, and set $\Meff_*:=\int a+Z\int b$.
\begin{enumerate}
\item[(i)] If $\Meff_*>0$, there is $C>0$ such that for all $\eps\in(0,1]$,
\[
T_\eps\le C\,\eps^{-\frac{2(p-1)}{2-N(p-1)}}\quad\text{if }p<\pF,\qquad
T_\eps\le \exp\big(C\eps^{-(p-1)}\big)\quad\text{if }p=\pF .
\]
\item[(ii)] If $\zeta$ satisfies \textup{(H$_a$)} (or, more generally, \textup{(H$_K$)} in \eqref{eq:HK}), there are $c,\eps_0>0$ such that for all $\eps\in(0,\eps_0]$,
\[
T_\eps\ge c\,\eps^{-\frac{2(p-1)}{2-N(p-1)}}\quad\text{if }p<\pF,\qquad
T_\eps\ge \exp\big(c\,\eps^{-(p-1)}\big)\quad\text{if }p=\pF ,
\]
without any sign condition on $a$, $b$ or $\Meff_*$.
\end{enumerate}
\end{theorem}

The upper bounds hold for weak solutions and use a test-function argument in the spirit of \cite{IkedaSobajima2019}; the lower bounds follow from the local version of Theorem~\ref{thm:exist}. The exponents coincide with those of Lee and Ni \cite{LeeNi1992} for the unforced problem with nonnegative data: integrable forcing does not alter the lifespan scale, only the constant through $\Meff_*$.

\subsection{Discussion and novelty}
The picture that emerges is the following. An integrable forcing never changes the critical exponent, and its whole effect in the Fujita range is encoded in the single number $\Meff$: a positive impulse $ZW$ subsidizes the initial mass, so that even data with $m_0<0$ blow up when $m_0>-ZW$ (Remark~\ref{rem:mass}). This contrasts sharply with the non-integrable regime, where the accumulated impulse diverges and shifts or destroys the threshold \cite{BLZ2000,JKS2020}. The mass-annihilation identity is a quantitative form of the test-function argument: rather than estimating the right-hand side of the weak formulation from above, we show that it tends to zero for every global solution, which is what allows us to reach the borderline case $\Meff=0$ and to describe the dynamics (vanishing of the $L^1$ norm) rather than merely nonexistence. In the unforced, nonpositive case it reduces to the classical extinction of mass for the absorption equation in the range $p\le\pF$, for which it provides a short proof (Remark~\ref{rem:absorption}).

On the existence side, the domination principle \eqref{eq:domination} is a simple but flexible device: it gives global existence for sign-changing data and forcing with a fully explicit smallness condition, it immediately provides the local existence time needed for sharp lifespan lower bounds, and it does not require integrability of $\zeta$ (Remark~\ref{rem:power}).

The paper is organized as follows. Section~\ref{sec:prelim} collects notation, the definitions of weak and mild solutions, and the auxiliary lemmas (heat kernel, cutoff functions, damped Duhamel formula, kernel estimates for the forcing). Section~\ref{sec:annihilation} proves Theorems~\ref{thm:annihilation} and \ref{thm:nonexist}. Section~\ref{sec:obstruction} is devoted to the localized obstruction and its consequences, including non-integrable forcings. Section~\ref{sec:existence} proves Theorems~\ref{thm:exist} and \ref{thm:nonpositive} and Corollary~\ref{cor:small}. Section~\ref{sec:lifespan} proves Theorem~\ref{thm:lifespan}, and Section~\ref{sec:asymptotics} proves Theorem~\ref{thm:asymptotics}. Section~\ref{sec:model} summarizes the complete picture for the model example \eqref{eq:model-zeta}, and Section~\ref{sec:open} lists open problems; all statements there are conjectural and are clearly separated from the results proved in the paper.

\section{Preliminaries}\label{sec:prelim}

\subsection{Notation and the heat semigroup}
For $1\le q\le\infty$ we write $\norm{\cdot}_q$ for the norm of $L^q(\RN)$, and $\tn{f}=\max\{\norm{f}_1,\norm{f}_\infty\}$. For $0<T\le\infty$ we set $\Q_T:=\RN\times(0,T)$, so that $\Q=\Q_\infty$. A function $u$ on $\Q_T$ belongs to $L^{q}_{\rm loc}(\RN\times[0,T))$ if $u\in L^{q}(K\cap\Q_T)$ for every compact set $K\subset\RN\times[0,T)$; such compact sets may touch the initial time. For a measurable $f$ on $\Q_T$ we use the Duhamel operator
\[
\Duh[f](x,t):=\int_0^{t}\!\!\int_{\RN}G_{t-s}(x-y)\,f(y,s)\,dy\ds,
\]
whenever the integral converges absolutely; we also write $\Duh[f](t)=\int_0^tS(t-s)f(s)\ds$. The cutoff exponent
\[
k:=2p'=\frac{2p}{p-1}>2
\]
will be used throughout. The letter $C$ denotes a positive constant whose dependence is indicated when relevant and which may change from line to line.

\begin{lemma}\label{lem:heat}
The following properties hold.
\begin{enumerate}
\item[(a)] $S(t)$ is positivity preserving, $\abs{S(t)f}\le S(t)\abs{f}$, $S(t)1=1$, and $S(t)S(s)=S(t+s)$ for $t,s\ge0$.
\item[(b)] $\norm{S(t)f}_q\le\norm{f}_q$ for $1\le q\le\infty$, $\norm{S(t)f}_\infty\le(4\pi t)^{-N/2}\norm{f}_1$; hence $\norm{S(t)f}_\infty\le\hN(t)\tn{f}$.
\item[(c)] $\int_{\RN}S(t)f\dx=\int_{\RN}f\dx$ for $f\in L^1(\RN)$.
\item[(d)] For $f\in L^1(\RN)$, $\norm{S(t)f-(\int f)\,G_t}_1\to0$ as $t\to\infty$.
\item[(e)] For every fixed $s\ge0$, $\norm{G_{t-s}-G_t}_1\to0$ as $t\to\infty$.
\item[(f)] $\hN$ is nonincreasing, $\hN(t/2)\le2^{N/2}\hN(t)$ for $t>0$, and $\int_0^\infty\hN(t)^{p-1}\dt<\infty$ if and only if $p>\pF$.
\end{enumerate}
\end{lemma}

\begin{proof}
Properties (a)--(c) are classical. For (d), assume first that $f\in C_c(\RN)$ with $\supp f\subset B_\rho$. Then
\[
\Big\|S(t)f-\Big(\int f\Big)G_t\Big\|_1\le\int_{B_\rho}\abs{f(y)}\,\norm{G_t(\cdot-y)-G_t}_1\,dy
=\int_{B_\rho}\abs{f(y)}\,\norm{G_1(\cdot-y/\sqrt t)-G_1}_1\,dy,
\]
which tends to zero by the continuity of translations in $L^1$, uniformly for $\abs{y}\le\rho$. For general $f\in L^1$ we approximate by $f_j\in C_c$ and use $\norm{S(t)(f-f_j)}_1+\abs{\int(f-f_j)}\norm{G_t}_1\le2\norm{f-f_j}_1$. For (e), the change of variables $x=\sqrt t\,y$ gives $\norm{G_{t-s}-G_t}_1=\norm{G_{(t-s)/t}-G_1}_1$, and $\lambda\mapsto G_\lambda$ is continuous from $(0,\infty)$ into $L^1$ by dominated convergence. In (f), the monotonicity is clear; if $2\pi t\ge1$ then $\hN(t/2)=(2\pi t)^{-N/2}=2^{N/2}\hN(t)$, while if $2\pi t<1$ then $\hN(t/2)=1<2^{N/2}\hN(t)$ because $(4\pi t)^{-N/2}>2^{-N/2}$. Finally $\hN^{p-1}$ is bounded and behaves like $t^{-N(p-1)/2}$ at infinity, which is integrable exactly when $N(p-1)/2>1$, that is $p>\pF$.
\end{proof}

\subsection{Weak and mild solutions}

\begin{definition}\label{def:weak}
Let $0<T\le\infty$, $u_0\in L^1_{\rm loc}(\RN)$, $w\in L^1_{\rm loc}(\RN)$ and $\zeta\in L^1_{\rm loc}([0,\infty))$. A function $u\in L^p_{\mathrm{loc}}(\RN\times[0,T))$ is called a \emph{weak solution} of \eqref{eq:model} on $\Q_T$ if, for every $\psi\in C_c^{\infty}\big(\RN\times[0,T)\big)$,
\begin{equation}\label{eq:weak}
\iint_{\Q_T} \abs{u}^{p}\,\psi \dx\dt
+\iint_{\Q_T} \zeta(t)\,w(x)\,\psi\dx\dt
+\int_{\RN}u_0(x)\,\psi(x,0)\dx
=-\iint_{\Q_T} u\,\big(\partial_t\psi+\Delta\psi\big)\dx\dt .
\end{equation}
It is a \emph{global weak solution} if $T=\infty$.
\end{definition}

All integrals in \eqref{eq:weak} are finite, since $\zeta\otimes w\in L^1_{\rm loc}(\RN\times[0,\infty))$ by Fubini's theorem. The test functions are compactly supported in space and time but \emph{need not vanish at} $t=0$; consequently the initial datum is genuinely felt through the boundary term in \eqref{eq:weak}. This is the point at which sign-changing data are correctly handled. We write $\Lop\psi:=\partial_t\psi+\Delta\psi$.

\begin{remark}\rm \label{rem:test}
Identity \eqref{eq:weak} extends to every $\psi\in C^{2,1}_c(\RN\times[0,T))$, that is, to compactly supported functions such that $\psi$, $\nabla\psi$, $D^2\psi$ and $\partial_t\psi$ are continuous on $\RN\times[0,T)$. Indeed, extend $\psi$ to negative times by the reflection $\psi(x,-t):=4\psi(x,t/2)-3\psi(x,t)$, $t>0$, which preserves the $C^{2,1}$ regularity across $t=0$, and mollify in $(x,t)$. The mollified functions $\psi_j$, restricted to $t\ge0$, belong to $C_c^\infty(\RN\times[0,T))$, have supports in a fixed compact subset of $\RN\times[0,T)$, and satisfy $\psi_j\to\psi$, $\Lop\psi_j\to\Lop\psi$ uniformly; since $u$, $\abs{u}^p$ and $\zeta\otimes w$ are integrable on that compact set and $u_0$ is integrable on its trace at $t=0$, we can pass to the limit in each term of \eqref{eq:weak}. All the cutoff functions used below are of class $C^{2,1}$, because they are powers $\eta^{k}$ with $k>2$ of a smooth function.
\end{remark}

\begin{definition}\label{def:mild}
Let $u_0,w\in L^1(\RN)\cap L^{\infty}(\RN)$, $\zeta\in L^1_{\rm loc}([0,\infty))$ and $0<T\le\infty$. A measurable function $u:\Q_T\to\R$ is a \emph{mild solution} of \eqref{eq:model} on $(0,T)$ if, for every $T'\in(0,T)$, $u\in L^\infty(\Q_{T'})$, $\operatorname*{ess\,sup}_{0<t<T'}\norm{u(t)}_1<\infty$, and
\begin{equation}\label{eq:mild}
u(t)=S(t)u_0+\int_0^{t}S(t-s)\abs{u(s)}^{p}\ds
+\int_0^{t}S(t-s)\big(\zeta(s)\,w\big)\ds
\end{equation}
for a.e. $(x,t)\in\Q_T$. It is a \emph{global mild solution} if $T=\infty$.
\end{definition}

The right-hand side of \eqref{eq:mild} does not depend on the representative of $u$ and is defined at every point of $\Q_T$; we always use it as the representative of $u$. With this convention $u(t)\in L^1\cap L^\infty$ for every $t\in(0,T)$, and the bounds in Definition~\ref{def:mild} hold for every $t$. Note also that $\abs{u}^p+\zeta w\in L^1(\Q_{T'})$, since $\norm{\abs{u(s)}^p}_1\le\norm{u}_{L^\infty(\Q_{T'})}^{p-1}\norm{u(s)}_1$.

\begin{lemma}\label{lem:mild-weak}
Let $u$ be a mild solution of \eqref{eq:model} on $(0,T)$ and set $f:=\abs{u}^p+\zeta w$.
\begin{enumerate}
\item[(a)] $u$ is a weak solution of \eqref{eq:model} on $\Q_T$.
\item[(b)] For every $t\in(0,T)$,
\begin{equation}\label{eq:mass-balance}
\int_{\RN}u(x,t)\dx=m_0+\int_0^{t}\!\!\int_{\RN}\abs{u}^{p}\dx\ds+W\int_0^t\zeta(s)\ds .
\end{equation}
\item[(c)] For $0<\tau<t<T$, $u(t)=S(t-\tau)u(\tau)+\int_\tau^tS(t-s)f(s)\ds$.
\end{enumerate}
\end{lemma}

\begin{proof}
Fix $\psi\in C_c^\infty(\RN\times[0,T))$ and $T'<T$ with $\supp\psi\subset\RN\times[0,T']$. We first claim that for every $g\in L^1(\RN)$ and $s\in[0,T)$,
\begin{equation}\label{eq:duality}
\int_s^{T}\!\!\int_{\RN}\big[S(t-s)g\big](x)\,\Lop\psi(x,t)\dx\dt=-\int_{\RN}g(x)\,\psi(x,s)\dx .
\end{equation}
By the symmetry of the heat kernel, $\int[S(t-s)g]\,\Lop\psi(t)\dx=\int g\,S(t-s)[\Lop\psi(t)]\dx$. For a function $\phi\in C^{2}$ with bounded derivatives one has $\frac{d}{d\tau}S(\tau)\phi=S(\tau)\Delta\phi$ uniformly, hence $\Theta(t):=\int g\,S(t-s)[\psi(t)]\dx$ is continuous on $[s,T)$, continuously differentiable on $(s,T)$ with $\Theta'(t)=\int g\,S(t-s)[\Lop\psi(t)]\dx$ (differentiation under the integral sign is justified by $g\in L^1$ and the uniform bounds on the derivatives of $\psi$), $\Theta(s)=\int g\,\psi(s)\dx$ because $S(\tau)\psi(s)\to\psi(s)$ uniformly as $\tau\to0$, and $\Theta(t)=0$ for $t\ge T'$. Integrating $\Theta'$ over $(s,T')$ gives \eqref{eq:duality}.

We now integrate \eqref{eq:mild} against $\Lop\psi$. Since $\norm{S(t-s)f(s)}_1\le\norm{f(s)}_1$ and $\int_0^{T'}\norm{f(s)}_1\ds<\infty$, Fubini's theorem applies and \eqref{eq:duality} yields
\[
\begin{aligned}
\iint_{\Q_T}u\,\Lop\psi
&=\int_0^{T}\!\!\int S(t)u_0\,\Lop\psi\dx\dt+\int_0^{T}\Big(\int_s^{T}\!\!\int S(t-s)f(s)\,\Lop\psi\dx\dt\Big)\ds\\
&=-\int u_0\,\psi(0)\dx-\iint_{\Q_T}f\,\psi ,
\end{aligned}
\]
which is \eqref{eq:weak}; moreover $u\in L^p_{\rm loc}(\RN\times[0,T))$ since $u$ is bounded on $\Q_{T'}$. This proves (a). Integrating \eqref{eq:mild} over $\RN$ and using Lemma~\ref{lem:heat}(c) and Fubini's theorem gives (b). For (c), the semigroup property and Fubini's theorem (all integrands are absolutely integrable) give $S(t-\tau)u(\tau)=S(t)u_0+\int_0^\tau S(t-s)f(s)\ds$, and it remains to add $\int_\tau^tS(t-s)f(s)\ds$.
\end{proof}

\begin{lemma}\label{lem:unique}
Two mild solutions of \eqref{eq:model} on $(0,T)$ coincide.
\end{lemma}

\begin{proof}
Let $u_1,u_2$ be mild solutions, fix $T'<T$, and let $K$ bound $\abs{u_1}$ and $\abs{u_2}$ on $\Q_{T'}$. Since $\abs{\abs{a}^p-\abs{b}^p}\le pK^{p-1}\abs{a-b}$ for $\abs{a},\abs{b}\le K$, Lemma~\ref{lem:heat}(a) gives $\abs{u_1-u_2}(t)\le pK^{p-1}\int_0^tS(t-s)\abs{u_1-u_2}(s)\ds$, hence $\phi(t):=\norm{(u_1-u_2)(t)}_\infty$ satisfies $\phi(t)\le pK^{p-1}\int_0^t\phi(s)\ds$ for $t<T'$. Iterating this inequality $n$ times yields $\phi(t)\le(pK^{p-1}t)^n\norm{\phi}_{L^\infty(0,T')}/n!$, and letting $n\to\infty$ gives $\phi\equiv0$ on $(0,T')$.
\end{proof}

The next lemma is the variation-of-constants formula with an additional linear damping, which turns monotone nonlinearities into order-preserving maps.

\begin{lemma}\label{lem:damped}
Let $0<T'<\infty$, $\lambda\ge0$, $v_0\in L^\infty(\RN)$, and let $F$ be measurable on $\Q_{T'}$ with $\int_0^{T'}\norm{F(s)}_\infty\ds<\infty$. For a bounded measurable function $v$ on $\Q_{T'}$, the identity
\begin{equation}\label{eq:duhamel-plain}
v(t)=S(t)v_0+\int_0^tS(t-s)F(s)\ds\qquad\text{a.e. in }\Q_{T'}
\end{equation}
holds if and only if
\begin{equation}\label{eq:duhamel-damped}
v(t)=e^{-\lambda t}S(t)v_0+\int_0^te^{-\lambda(t-s)}S(t-s)\big(F(s)+\lambda v(s)\big)\ds\qquad\text{a.e. in }\Q_{T'} .
\end{equation}
\end{lemma}

\begin{proof}
Assume \eqref{eq:duhamel-plain}. Inserting it in the term containing $\lambda v$, using the semigroup property and Fubini's theorem (justified by $\int_0^{T'}\norm{F(s)}_\infty\ds<\infty$), we find
\[
\lambda\int_0^te^{-\lambda(t-s)}S(t-s)v(s)\ds
=\big(1-e^{-\lambda t}\big)S(t)v_0+\int_0^t\big(1-e^{-\lambda(t-r)}\big)S(t-r)F(r)\,dr,
\]
because $\lambda\int_r^te^{-\lambda(t-s)}\ds=1-e^{-\lambda(t-r)}$. Adding $e^{-\lambda t}S(t)v_0+\int_0^te^{-\lambda(t-s)}S(t-s)F(s)\ds$ to both sides gives \eqref{eq:duhamel-damped}. Conversely, assume \eqref{eq:duhamel-damped} and let $\tilde v$ be the right-hand side of \eqref{eq:duhamel-plain}, which is bounded on $\Q_{T'}$ by $\norm{v_0}_\infty+\int_0^{T'}\norm{F}_\infty$. By the first part, $\tilde v$ satisfies \eqref{eq:duhamel-damped}, so that $z:=v-\tilde v$ satisfies $z(t)=\lambda\int_0^te^{-\lambda(t-s)}S(t-s)z(s)\ds$ and therefore $\norm{z(t)}_\infty\le\lambda\int_0^t\norm{z(s)}_\infty\ds$. As in the proof of Lemma~\ref{lem:unique}, iteration gives $z=0$.
\end{proof}

\subsection{Cutoff functions}
We fix once and for all a function
\begin{equation}\label{eq:eta}
\eta\in C^\infty([0,\infty)),\qquad 0\le\eta\le1,\qquad \eta\equiv1\text{ on }[0,1],\qquad \eta\equiv0\text{ on }[2,\infty),\qquad \eta'\le0 .
\end{equation}
We shall use two families of test functions built on $\eta$. The \emph{separated cutoffs}, centred at $x_0\in\RN$, are
\begin{equation}\label{eq:Psi}
\Psi_{R,T}(x,t):=\eta\!\Big(\tfrac{\abs{x-x_0}^2}{R^2}\Big)^{k}\,
\eta\!\Big(\tfrac{t}{T}\Big)^{k}\in C^{2,1}_c\big(\RN\times[0,\infty)\big),
\qquad \Psi_{R,T}(x,0)=\eta\!\Big(\tfrac{\abs{x-x_0}^2}{R^2}\Big)^{k},
\end{equation}
and the \emph{parabolic cutoffs} are
\begin{equation}\label{eq:Psi-par}
\Psi_R(x,t):=\eta\big(s_R(x,t)\big)^{k},\qquad s_R(x,t):=\frac{\abs{x}^2+t}{R^2},\qquad \Psi_R(x,0)=\eta\Big(\tfrac{\abs x^2}{R^2}\Big)^k .
\end{equation}
The choice $\eta(0)=1$ is essential: it makes the boundary term in \eqref{eq:weak} nontrivial and lets the forcing contribute its full impulse. In both lemmas below the integrands $\abs{\Lop\Psi}^{p'}\Psi^{1-p'}$ are understood to vanish where $\Psi=0$; this is consistent, since $\Lop\Psi$ contains the factor $\eta^{k-2}$ and therefore vanishes wherever $\Psi$ does.

\begin{lemma}\label{lem:cutoff}
There is a constant $C=C(N,p,\eta)>0$ such that for all $x_0\in\RN$ and $R,T>0$,
\begin{align}
\iQ \frac{\abs{\partial_t\Psi_{R,T}}^{p'}}{\Psi_{R,T}^{\,p'-1}}\dx\dt
&\le C\,R^{N}\,T^{1-p'},\label{eq:time-cost}\\[2pt]
\iQ \frac{\abs{\Delta\Psi_{R,T}}^{p'}}{\Psi_{R,T}^{\,p'-1}}\dx\dt
&\le C\,R^{N-2p'}\,T.\label{eq:space-cost}
\end{align}
Consequently, with $|a+b|^{p'}\le 2^{p'-1}(|a|^{p'}+|b|^{p'})$,
\begin{equation}\label{eq:IRT}
I_{R,T}:=\iQ\frac{\abs{\partial_t\Psi_{R,T}+\Delta\Psi_{R,T}}^{p'}}{\Psi_{R,T}^{\,p'-1}}\dx\dt
\le C\big(R^{N}T^{1-p'}+R^{N-2p'}T\big),
\end{equation}
and the parabolic balance $T=R^2$ yields $I_{R,R^2}\le C\,R^{\,N+2-2p'}=C\,R^{N-\frac{2}{p-1}}$.
\end{lemma}

\begin{proof}
Write $\phi_R(x):=\eta(\abs{x-x_0}^2/R^2)$ and $\eta_T(t):=\eta(t/T)$, so that $\Psi_{R,T}=\phi_R^k\eta_T^k$. Since $\partial_t\Psi_{R,T}=k\,\phi_R^{k}\eta_T^{k-1}\eta_T'$ and $k=2p'$,
\[
\begin{aligned}
\frac{\abs{\partial_t\Psi_{R,T}}^{p'}}{\Psi_{R,T}^{\,p'-1}}
&=k^{p'}\,\phi_R^{\,k}\,\eta_T^{\,(k-1)p'-k(p'-1)}\,\abs{\eta_T'}^{p'}
=k^{p'}\,\phi_R^{\,k}\,\eta_T^{\,p'}\,\abs{\eta_T'}^{p'}\\
&\le k^{p'}\norm{\eta'}_\infty^{p'}\,T^{-p'}\,\mathbf 1_{\{\abs{x-x_0}\le\sqrt2R\}}\mathbf 1_{\{T\le t\le2T\}},
\end{aligned}
\]
because $(k-1)p'-k(p'-1)=k-p'=p'$ and $\eta_T'(t)=T^{-1}\eta'(t/T)$ vanishes outside $[T,2T]$. Integrating over a ball of radius $\sqrt2R$ and an interval of length $T$ proves \eqref{eq:time-cost}.

For the Laplacian, $\Delta\Psi_{R,T}=\eta_T^{k}\,\Delta(\phi_R^{k})$ and
\[
\Delta(\phi_R^{k})=k\,\phi_R^{k-1}\Delta\phi_R+k(k-1)\,\phi_R^{k-2}\abs{\nabla\phi_R}^2 .
\]
With $r:=\abs{x-x_0}^2/R^2$ we have $\nabla\phi_R=2(x-x_0)R^{-2}\eta'(r)$ and $\Delta\phi_R=R^{-2}\big(2N\eta'(r)+4r\eta''(r)\big)$; both vanish unless $1\le r\le2$, so that $\abs{\nabla\phi_R}^2\le 8\norm{\eta'}_\infty^2R^{-2}$ and $\abs{\Delta\phi_R}\le(2N\norm{\eta'}_\infty+8\norm{\eta''}_\infty)R^{-2}$. Using $\phi_R\le1$ we obtain $\abs{\Delta(\phi_R^{k})}\le C R^{-2}\phi_R^{k-2}\mathbf 1_{\{R\le\abs{x-x_0}\le\sqrt2 R\}}$, and after the exact cancellation $(k-2)p'-k(p'-1)=k-2p'=0$,
\[
\frac{\abs{\Delta\Psi_{R,T}}^{p'}}{\Psi_{R,T}^{\,p'-1}}
=\eta_T^{\,k}\,\frac{\abs{\Delta(\phi_R^{k})}^{p'}}{\phi_R^{\,k(p'-1)}}
\le C\,R^{-2p'}\,\mathbf 1_{\{R\le\abs{x-x_0}\le\sqrt2R\}}\mathbf 1_{\{0\le t\le 2T\}} .
\]
Integrating gives \eqref{eq:space-cost}, and \eqref{eq:IRT} follows. Finally $N+2-2p'=N-2/(p-1)$ because $2p'-2=2/(p-1)$.
\end{proof}

\begin{lemma}\label{lem:cutoff-par}
For every $R>0$ the function $\Psi_R$ belongs to $C^{2,1}_c(\RN\times[0,\infty))$, is supported in $\{\abs{x}^2+t\le2R^2\}$, is nondecreasing in $R$, and $\Psi_R\to1$ pointwise as $R\to\infty$. Moreover there is $C_1=C_1(N,p,\eta)$ such that
\begin{equation}\label{eq:LPsi-par}
\abs{\Lop\Psi_R}\le C_1\,R^{-2}\,\eta(s_R)^{k-2}\,\mathbf 1_{\{1\le s_R\le2\}} .
\end{equation}
Consequently, for every $u\in L^p_{\rm loc}(\RN\times[0,\infty))$ and $R>0$,
\begin{equation}\label{eq:Holder-par}
\iQ\abs{u}\,\abs{\Lop\Psi_R}\dx\dt\le C_2\,R^{\frac{N+2-2p'}{p'}}\Big(\iint_{\mathcal D_R}\abs{u}^{p}\,\Psi_R\dx\dt\Big)^{1/p},
\end{equation}
where $\mathcal D_R:=\{(x,t)\in\Q:\ R^2\le\abs{x}^2+t\le2R^2\}$,
with $C_2=C_2(N,p,\eta)$, and
\begin{equation}\label{eq:cost-par}
\iQ\frac{\abs{\Lop\Psi_R}^{p'}}{\Psi_R^{\,p'-1}}\dx\dt\le C_2^{p'}\,R^{N+2-2p'} .
\end{equation}
\end{lemma}

\begin{proof}
The monotonicity in $R$ follows from $\eta'\le0$, and $\Psi_R\to\eta(0)^k=1$. Put $g:=\eta^k$, so that $\Psi_R=g(s_R)$. Since $\partial_ts_R=R^{-2}$, $\nabla s_R=2xR^{-2}$ and $\Delta s_R=2NR^{-2}$,
\[
\Lop\Psi_R=R^{-2}\Big[(2N+1)\,g'(s_R)+4\,\frac{\abs{x}^2}{R^2}\,g''(s_R)\Big].
\]
Now $g'=k\eta^{k-1}\eta'$ and $g''=k(k-1)\eta^{k-2}(\eta')^2+k\eta^{k-1}\eta''$ vanish outside $[1,2]$ and are bounded by a constant times $\eta^{k-2}$ there; since $\abs{x}^2/R^2\le s_R\le2$ on the support, \eqref{eq:LPsi-par} follows. On $\{1\le s_R\le2\}$ we write $\eta^{k-2}=(\eta^{k})^{1/p}\eta^{k/p'-2}=\Psi_R^{1/p}$, since $k/p'=2$. Hence, by H\"older's inequality,
\[
\iQ\abs{u}\abs{\Lop\Psi_R}\le C_1R^{-2}\iint_{\mathcal D_R}\abs{u}\Psi_R^{1/p}
\le C_1R^{-2}\Big(\iint_{\mathcal D_R}\abs{u}^p\Psi_R\Big)^{1/p}\abs{\mathcal D_R}^{1/p'} ,
\]
and $\abs{\mathcal D_R}\le\abs{B_1}(\sqrt2R)^N\cdot2R^2$. This gives \eqref{eq:Holder-par} because $-2+(N+2)/p'=(N+2-2p')/p'$. Similarly, $\abs{\Lop\Psi_R}^{p'}\Psi_R^{1-p'}\le C_1^{p'}R^{-2p'}\eta(s_R)^{(k-2)p'-k(p'-1)}\mathbf 1_{\mathcal D_R}=C_1^{p'}R^{-2p'}\mathbf 1_{\mathcal D_R}$, which gives \eqref{eq:cost-par}.
\end{proof}

\subsection{The forcing kernel}
The behavior of the linear forced flow is governed by the function $\abs{\zeta}\conv\hN$ defined in \eqref{eq:hN}. Besides (H$_a$) we shall use the following weaker hypothesis:
\begin{equation}\label{eq:HK}
\text{(H$_K$)}\qquad (\abs{\zeta}\conv\hN)(t)\le K\,\hN(t)\quad\text{for all }t>0,\ \text{for some constant }K\ge0 .
\end{equation}

\begin{lemma}\label{lem:kernel}
~\begin{enumerate}
\item[(i)] If $\zeta$ satisfies \textup{(H$_a$)}, then $\zeta$ satisfies \textup{(H$_K$)} and
\begin{equation}\label{eq:o-cond}
\lim_{t\to\infty}t^{N/2}\int_{t/2}^{t}\abs{\zeta(s)}\,\hN(t-s)\ds=0 .
\end{equation}
\item[(ii)] If $N\ge3$, $\zeta\in L^1(0,\infty)$ and $\abs{\zeta(t)}\le Ct^{-b}$ for a.e.\ $t\ge1$ with $b>1$, then $(\abs{\zeta}\conv\hN)(t)\le C'(1+t)^{-\min\{b,N/2\}}$ for all $t>0$.
\item[(iii)] If $q\ge1$ and $\zeta\in L^1(0,\infty)$, then $\norm{\abs{\zeta}\conv\hN}_{L^q(0,\infty)}\le\norm{\zeta}_{L^1}\norm{\hN}_{L^q(0,\infty)}$.
\end{enumerate}
\end{lemma}

\begin{proof}
(i) For $t\ge2$ we split the integral defining $(\abs{\zeta}\conv\hN)(t)$ at $s=t/2$. On $[0,t/2]$ we have $t-s\ge t/2$, so that by Lemma~\ref{lem:heat}(f)
\[
\int_0^{t/2}\abs{\zeta(s)}\hN(t-s)\ds\le\hN(t/2)\norm{\zeta}_1\le2^{N/2}\norm{\zeta}_1\hN(t).
\]
On $[t/2,t]$ we have $s\ge1$ and $\abs{\zeta(s)}\le C_\zeta2^at^{-a}$, hence
\[
\int_{t/2}^{t}\abs{\zeta(s)}\hN(t-s)\ds\le C_\zeta2^a\,t^{-a}\int_0^{t/2}\hN(r)\,dr .
\]
The last integral is bounded by $1+(t/2\pi)^{1/2}$ if $N=1$, by $C(1+\log(1+t))$ if $N=2$, and by a constant if $N\ge3$. Since $a>1$ when $N\le2$ and $a>N/2$ when $N\ge3$, we obtain \eqref{eq:o-cond}. In particular $t^{N/2}(\abs{\zeta}\conv\hN)(t)$ is bounded on $[2,\infty)$, and since $\hN(t)=(4\pi t)^{-N/2}$ there, $(\abs{\zeta}\conv\hN)\le K_1\hN$ on $[2,\infty)$. On $(0,2]$ we simply use $(\abs{\zeta}\conv\hN)(t)\le\norm{\zeta}_1\le\norm{\zeta}_1\hN(2)^{-1}\hN(t)$.

(ii) The same splitting gives, for $t\ge2$, $(\abs{\zeta}\conv\hN)(t)\le2^{N/2}\norm{\zeta}_1\hN(t)+C2^bt^{-b}\norm{\hN}_{L^1(0,\infty)}$, and $\hN\in L^1(0,\infty)$ because $N\ge3$; for $t\le2$ the bound $\norm{\zeta}_1$ suffices.

(iii) This is Young's inequality on the real line, applied to $\abs{\zeta}\mathbf 1_{(0,\infty)}$ and $\hN\mathbf 1_{(0,\infty)}$, whose convolution coincides with $\abs{\zeta}\conv\hN$ on $(0,\infty)$.
\end{proof}

\section{Mass annihilation in the Fujita range}\label{sec:annihilation}

\subsection{Proof of Theorem~\ref{thm:annihilation}}
Let $u$ be a global weak solution and fix $R\ge1$. By Remark~\ref{rem:test} and Lemma~\ref{lem:cutoff-par} we may test \eqref{eq:weak} with the parabolic cutoff $\psi=\Psi_R$ of \eqref{eq:Psi-par}. Since $\Psi_R(x,0)=\eta(\abs{x}^2/R^2)^{k}$, the identity reads
\begin{equation}\label{eq:identity}
\underbrace{\iQ\abs{u}^{p}\Psi_{R}}_{=:\,J(R)}
+\underbrace{\iQ \zeta\,w\,\Psi_{R}}_{=:\,\mathcal F_{R}}
+\underbrace{\int_{\RN}u_0\,\eta\big(\abs{x}^2/R^2\big)^{k}\dx}_{=:\,\mathcal I_{R}}
=-\iQ u\,\Lop\Psi_{R} .
\end{equation}
All four quantities are finite because $\Psi_R$ has compact support.

\emph{The data terms.} Since $0\le\Psi_R\le1$, $\Psi_R\to1$ pointwise, $u_0\in L^1(\RN)$ and $\zeta\otimes w\in L^1(\Q)$ (with $\norm{\zeta\otimes w}_{L^1(\Q)}=\norm{\zeta}_1\norm{w}_1$), dominated convergence gives
\begin{equation}\label{eq:limits}
\mathcal I_R\xrightarrow[R\to\infty]{}m_0,
\qquad
\mathcal F_R\xrightarrow[R\to\infty]{}\iQ\zeta(t)w(x)\dx\dt=Z\,W,
\end{equation}
and therefore
\begin{equation}\label{eq:F+I-limit}
\mathcal F_{R}+\mathcal I_R\xrightarrow[R\to\infty]{}m_0+Z\,W=\Meff,
\qquad
\abs{\mathcal F_{R}+\mathcal I_R}\le K_0:=\norm{u_0}_1+\norm{\zeta}_1\norm{w}_1 .
\end{equation}

\emph{A quantitative inequality.} By \eqref{eq:Holder-par} and Young's inequality, for every $R>0$,
\[
\Big|\iQ u\,\Lop\Psi_R\Big|\le C_2R^{\frac{N+2-2p'}{p'}}J(R)^{1/p}\le\frac1p\,J(R)+\frac{C_2^{p'}}{p'}\,R^{N+2-2p'} .
\]
Inserting this into \eqref{eq:identity} and using $1-\frac1p=\frac1{p'}$ we obtain, for every $p>1$ and every $R>0$,
\begin{equation}\label{eq:subcritical}
\frac1{p'}\,J(R)+\mathcal F_R+\mathcal I_R\le\frac{C_2^{p'}}{p'}\,R^{\,N+2-2p'} .
\end{equation}
When $1<p<\pF$ we have $N+2-2p'<0$, and letting $R\to\infty$ in \eqref{eq:subcritical} already gives $\Meff\le0$; to obtain the identity \eqref{eq:annihilation}, and to include $p=\pF$, we argue as follows.

\emph{Integrability of $\abs{u}^{p}$.} Since $p\le\pF$ means $N+2-2p'\le0$, we have $R^{N+2-2p'}\le1$ for $R\ge1$, and \eqref{eq:subcritical} together with \eqref{eq:F+I-limit} gives $J(R)\le p'K_0+C_2^{p'}=:K_1$ for all $R\ge1$. As $R\nearrow\infty$ the cutoffs $\Psi_R$ increase to $1$ (Lemma~\ref{lem:cutoff-par}), so monotone convergence yields
\begin{equation}\label{eq:LpQ}
\iQ\abs{u}^{p}\dx\dt\le K_1<\infty ,\qquad J(R)\xrightarrow[R\to\infty]{}\iQ\abs{u}^p\dx\dt .
\end{equation}

\emph{Vanishing of the diffusive term.} The function $\Lop\Psi_R$ is supported in the annular region $\mathcal D_R=\{R^2\le\abs{x}^2+t\le2R^2\}$, which escapes to infinity as $R\to\infty$. By \eqref{eq:LpQ} and dominated convergence,
\begin{equation}\label{eq:tail}
\iint_{\mathcal D_{R}}\abs{u}^{p}\Psi_R\dx\dt
\le\iint_{\{\abs{x}^2+t\ge R^2\}}\abs{u}^{p}\dx\dt\xrightarrow[R\to\infty]{}0 .
\end{equation}
Using \eqref{eq:Holder-par} once more, with $R^{(N+2-2p')/p'}\le1$,
\[
\Big|\iQ u\,\Lop\Psi_R\Big|\le C_2\Big(\iint_{\mathcal D_{R}}\abs{u}^{p}\Psi_R\Big)^{1/p}\xrightarrow[R\to\infty]{}0 .
\]

\emph{Conclusion.} Letting $R\to\infty$ in \eqref{eq:identity} and using \eqref{eq:F+I-limit}, \eqref{eq:LpQ} and the last display, we obtain $\iQ\abs{u}^p+\Meff=0$, which is \eqref{eq:annihilation}. \qed

\subsection{Proof of Theorem~\ref{thm:nonexist}}
Statements (i) and (iii) are immediate consequences of \eqref{eq:annihilation}, whose left-hand side is nonnegative. Assume $\Meff=0$ and let $u$ be a global weak solution. Then \eqref{eq:annihilation} gives $u=0$ a.e. in $\Q$, and \eqref{eq:weak} reduces to
\[
\iQ\zeta(t)w(x)\psi\dx\dt+\int_{\RN}u_0\,\psi(x,0)\dx=0\qquad\text{for all }\psi\in C_c^\infty(\RN\times[0,\infty)).
\]
Choosing first $\psi\in C^\infty_c(\Q)$, which vanishes near $t=0$, we see that the locally integrable function $\zeta\otimes w$ vanishes as a distribution in $\Q$, hence a.e. Choosing then $\psi(x,t)=\varphi(x)\chi(t)$ with $\varphi\in C^\infty_c(\RN)$ and $\chi\in C_c^\infty([0,\infty))$, $\chi(0)=1$, we get $\int u_0\varphi=0$ for all $\varphi$, that is $u_0=0$ a.e. Conversely, if $u_0=0$ and $\zeta\otimes w=0$ a.e., then $u\equiv0$ is a global weak solution, and any global weak solution satisfies $\iQ\abs{u}^p=-\Meff=0$ by \eqref{eq:annihilation}, so that $u\equiv0$ is the only one. By Fubini's theorem, $\zeta(t)w(x)=0$ a.e. in $\Q$ if and only if $\zeta=0$ a.e. in $(0,\infty)$ or $w=0$ a.e. in $\RN$. \qed

\subsection{Comments}

\begin{remark}\rm \label{rem:mass}
Theorems~\ref{thm:annihilation} and \ref{thm:nonexist} isolate exactly what an integrable forcing does. The critical exponent is the Fujita exponent $\pF$ of the unforced source equation $\partial_t u-\Delta u=\abs{u}^{p}$; the forcing does \emph{not} change it. What the forcing changes is the admissibility threshold on the data: in place of the bare condition $m_0>0$ one obtains the \emph{effective mass} condition \eqref{eq:mass-condition}, $m_0+Z\,W>0$, and its borderline version $m_0+ZW=0$. A positive impulse $ZW>0$ subsidizes the initial mass, so that even sign-changing data with $m_0\le0$ must blow up provided $m_0\ge-ZW$ and the data are nontrivial; conversely a negative impulse $ZW<0$ can be compensated by a positive initial mass. In the model case \eqref{eq:model-zeta} the parameter $\sigma$ enters only through the constant $Z=\Gamma(\sigma+1)$. Setting $w\equiv0$ recovers the Fujita statement $m_0\ge0$, $u_0\not\equiv0$ $\Rightarrow$ no global solution for $1<p\le\pF$.
\end{remark}

\begin{remark}\rm \label{rem:why-pF}
The proof of Theorem~\ref{thm:annihilation} uses $p\le\pF$ only through the boundedness of the diffusive cost $R^{N+2-2p'}$ in \eqref{eq:subcritical} and \eqref{eq:Holder-par}. For $p>\pF$ the identity \eqref{eq:annihilation} is false in general: by Corollary~\ref{cor:small} and Theorem~\ref{thm:asymptotics}(iii), small global solutions with $\Meff>0$ exist and have asymptotic mass $m_\infty=\Meff+\iQ\abs{u}^p>0$.
\end{remark}

\begin{remark}\rm \label{rem:absorption}
If $u_0\le0$ and $\zeta w\le0$, Theorem~\ref{thm:nonpositive} provides a global solution $u\le0$, and $v:=-u\ge0$ solves the absorption problem
\[
\partial_tv-\Delta v+v^{p}=\abs{\zeta(t)w(x)},\qquad v(0)=\abs{u_0}.
\]
For $1<p\le\pF$, Theorem~\ref{thm:annihilation} then states that $\iQ v^{p}=\norm{v_0}_1+\norm{\zeta w}_{L^1(\Q)}$: all the mass that is present initially or injected later is eventually absorbed, and by Theorem~\ref{thm:asymptotics}(ii) $\norm{v(t)}_1\to0$. For $\zeta w\equiv0$ this is the classical extinction of mass for the heat equation with absorption in the range $p\le\pF$ \cite{GmiraVeron1984,KaminPeletier1985,GalaktionovKurdyumovSamarskii1986,EscobedoKavianMatano1995}, for which the present argument gives a short proof based only on the weak formulation.
\end{remark}

\section{A localized obstruction for every \texorpdfstring{$p>1$}{p>1}}\label{sec:obstruction}

\subsection{Proof of Theorem~\ref{thm:obstruction}}
Let $u$ be a global weak solution, and let $x_0\in\RN$, $R,T>0$. We test \eqref{eq:weak} with the separated cutoff $\Psi_{R,T}$ of \eqref{eq:Psi}, which is admissible by Remark~\ref{rem:test}. Writing $\phi_R(x)=\eta(\abs{x-x_0}^2/R^2)$ and $\eta_T(t)=\eta(t/T)$, we obtain
\begin{equation}\label{eq:identity-sep}
J_{R,T}+\mathcal F_{R,T}+\mathcal I_{R}=-\iQ u\,\Lop\Psi_{R,T},
\qquad J_{R,T}:=\iQ\abs{u}^p\Psi_{R,T},\quad \mathcal I_R:=\int_{\RN}u_0\,\phi_R^{k}\dx ,
\end{equation}
and, because $\Psi_{R,T}$ is a product, the forcing term separates:
\begin{equation}\label{eq:forcing-split}
\mathcal F_{R,T}:=\iQ\zeta\,w\,\Psi_{R,T}
=\Big(\int_0^{\infty}\zeta(t)\,\eta_T(t)^{k}\dt\Big)
\Big(\int_{\RN}w(x)\,\phi_R(x)^{k}\dx\Big).
\end{equation}
Since $\Lop\Psi_{R,T}$ vanishes where $\Psi_{R,T}$ does, H\"older's inequality (in the weight $\Psi_{R,T}$) and Young's inequality give
\[
\Big|\iQ u\,\Lop\Psi_{R,T}\Big|
\le J_{R,T}^{1/p}\,I_{R,T}^{1/p'}
\le\frac1p\,J_{R,T}+\frac1{p'}\,I_{R,T},
\]
with $I_{R,T}$ as in \eqref{eq:IRT}. Inserting this into \eqref{eq:identity-sep} gives $\frac1{p'}J_{R,T}+\mathcal F_{R,T}+\mathcal I_R\le\frac1{p'}I_{R,T}$; dropping the nonnegative first term and invoking Lemma~\ref{lem:cutoff} we obtain \eqref{eq:obstruction} with $C_0:=C/p'$, where $C$ is the constant of \eqref{eq:IRT}. \qed

\medskip
With the parabolic choice $T=R^2$, \eqref{eq:obstruction} becomes
\begin{equation}\label{eq:obstruction-par}
\int_{\RN}u_0\,\phi_R^{k}\dx+\Big(\int_0^{\infty}\zeta\,\eta\big(\tfrac{t}{R^2}\big)^{k}\dt\Big)\Big(\int_{\RN}w\,\phi_R^{k}\dx\Big)\le 2C_0\,R^{\,N-\frac{2}{p-1}} ,
\end{equation}
and we note the trivial bounds $\abs{\int_0^\infty\zeta\eta_T^{k}}\le\norm{\zeta}_{L^1(0,2T)}$ and $\abs{\int w\phi_R^{k}}\le\norm{w}_{L^1(B(x_0,\sqrt2R))}$.

\subsection{Consequences}

\begin{corollary}\label{cor:obstruction}
Let $p>1$.
\begin{enumerate}
\item[(a)] \textup{(Large data.)} Let $\zeta\in L^1_{\rm loc}([0,\infty))$ and $w\in L^1_{\rm loc}(\RN)$ be fixed, and let $g\in L^1_{\rm loc}(\RN)$ satisfy $\abs{\{g>0\}}>0$. There is $\lambda_*>0$ such that \eqref{eq:model} with $u_0=\lambda g$ has no global weak solution for $\lambda>\lambda_*$.
\item[(b)] \textup{(Large forcing.)} Let $u_0\in L^1_{\rm loc}(\RN)$ be fixed, let $\zeta\in L^1_{\rm loc}([0,\infty))$ be nonnegative and not a.e. zero, and let $w\in L^1_{\rm loc}(\RN)$ satisfy $\abs{\{w>0\}}>0$. There is $\mu_*>0$ such that \eqref{eq:model} with forcing $\mu\zeta(t)w(x)$ has no global weak solution for $\mu>\mu_*$.
\item[(c)] \textup{(Slowly decaying data.)} Let $\zeta\in L^1(0,\infty)$, $w\in L^1(\RN)$, and let $u_0\in L^1_{\rm loc}(\RN)$ satisfy $u_0(x)\ge A\abs{x}^{-\gamma}$ for a.e. $\abs{x}\ge R_0$, for some $A,R_0>0$ and $\gamma\ge0$. If $\gamma<\min\{N,2/(p-1)\}$, there is no global weak solution. If $\gamma=2/(p-1)<N$, there is $A_*=A_*(N,p,\eta)$ such that there is no global weak solution when $A>A_*$.
\end{enumerate}
\end{corollary}

\begin{proof}
(a) Almost every point of $\{g>0\}$ is a Lebesgue point of $g$; fix such a point $x_0$. Then
\[
\Big|R^{-N}\!\int_{\RN}\big(g(x)-g(x_0)\big)\phi_R(x)^k\dx\Big|\le R^{-N}\!\int_{B(x_0,\sqrt2R)}\abs{g-g(x_0)}\dx\xrightarrow[R\to0]{}0,
\]
while $R^{-N}\int\phi_R^k\dx=c_\eta:=\int\eta(\abs{y}^2)^kdy>0$. Hence there is $R_0>0$ with $G_0:=\int g\,\phi_{R_0}^k\dx\ge\frac12c_\eta\,g(x_0)R_0^N>0$. If a global weak solution with $u_0=\lambda g$ exists, \eqref{eq:obstruction-par} with $R=R_0$ gives
\[
\lambda G_0\le2C_0R_0^{\,N-\frac2{p-1}}+\norm{\zeta}_{L^1(0,2R_0^2)}\norm{w}_{L^1(B(x_0,\sqrt2R_0))}=:\lambda_*G_0 .
\]

(b) As in (a), choose $x_0,R_0$ with $w_0:=\int w\,\phi_{R_0}^k\dx>0$, and choose $T_0>0$ with $\int_0^{T_0}\zeta>0$, so that $z_0:=\int_0^\infty\zeta\eta_{T_0}^k\ge\int_0^{T_0}\zeta>0$. If a global weak solution exists with forcing $\mu\zeta w$, \eqref{eq:obstruction} gives $\mu z_0w_0\le C_0(R_0^NT_0^{1-p'}+R_0^{N-2p'}T_0)+\norm{u_0}_{L^1(B(x_0,\sqrt2R_0))}$, which bounds $\mu$.

(c) Take $x_0=0$ and $R\ge R_0$ in \eqref{eq:obstruction-par}. Since $\phi_R=1$ on $B_R$ and $u_0\ge0$ outside $B_{R_0}$,
\[
\int u_0\phi_R^k\dx\ge-\norm{u_0}_{L^1(B_{R_0})}+A\int_{R_0\le\abs{x}\le R}\abs{x}^{-\gamma}\dx
=-\norm{u_0}_{L^1(B_{R_0})}+\frac{A\abs{\mathbb S^{N-1}}}{N-\gamma}\big(R^{N-\gamma}-R_0^{N-\gamma}\big),
\]
while the forcing term is bounded below by $-\norm{\zeta}_1\norm{w}_1$. Dividing \eqref{eq:obstruction-par} by $R^{N-\gamma}$ and letting $R\to\infty$, we get $A\abs{\mathbb S^{N-1}}/(N-\gamma)\le\limsup_{R\to\infty}2C_0R^{\gamma-2/(p-1)}$. The right-hand side vanishes if $\gamma<2/(p-1)$, a contradiction; if $\gamma=2/(p-1)$ it equals $2C_0$, and we reach a contradiction when $A>A_*:=2C_0(N-\gamma)/\abs{\mathbb S^{N-1}}$.
\end{proof}

Part (c) is the forced counterpart of the result of Lee and Ni \cite{LeeNi1992} for slowly decaying data; an integrable forcing does not affect it. In the Fujita range $p\le\pF$, the condition in (c) reads $\gamma<N$, so that $u_0\notin L^1$ and, heuristically, $\Meff=+\infty$, in accordance with Theorem~\ref{thm:nonexist}.

\begin{example}\label{ex:both}
Let $1<p\le\pF$, let $a\in C_c(B_1)$ be nonnegative and not identically zero, fix $x_1\in\RN$ with $\abs{x_1}\ge4$, and consider $\zeta\in L^1(0,\infty)$ and $w\in L^1(\RN)\cap L^\infty(\RN)$ with $\zeta\ge0$, $w\le0$. For $\lambda>0$ the data $u_0=-\lambda a$ have $\Meff=-\lambda\int a+ZW<0$, and Theorem~\ref{thm:nonpositive} gives a global solution. The data $u_0=\lambda\big(a-2a(\cdot-x_1)\big)$ also have $\Meff=-\lambda\int a+ZW<0$, but $\{u_0>0\}$ has positive measure, and Corollary~\ref{cor:obstruction}(a) shows that there is no global weak solution for $\lambda$ large.
\end{example}

\subsection{Non-integrable forcing}
The obstruction \eqref{eq:obstruction} also explains why the integrability of $\zeta$ is the dividing line between the present theory and the results for power forcings.

\begin{corollary}\label{cor:nonintegrable}
Let $u_0\in L^1(\RN)$, $w\in L^1(\RN)$ with $W>0$, and $\zeta\in L^1_{\rm loc}([0,\infty))$ with negative part $\zeta_-\in L^1(0,\infty)$; set $Z(T):=\int_0^T\zeta\dt$.
\begin{enumerate}
\item[(i)] If $1<p\le\pF$ and $Z(T)\to+\infty$ as $T\to\infty$, then \eqref{eq:model} has no global weak solution, whatever the size of $u_0$ and $w$.
\item[(ii)] If $p>\pF$ and $\limsup_{T\to\infty}T^{\frac1{p-1}-\frac N2}Z(T)=+\infty$, then \eqref{eq:model} has no global weak solution.
\end{enumerate}
\end{corollary}

\begin{proof}
Take $x_0=0$ in \eqref{eq:obstruction-par}. Writing $\zeta=\zeta_+-\zeta_-$ and using $\eta_T^k=1$ on $[0,T]$,
\[
\int_0^\infty\zeta\,\eta_T^{k}\dt\ge\int_0^T\zeta_+\dt-\norm{\zeta_-}_1\ge Z(T)-\norm{\zeta_-}_1 .
\]
Since $\int w\phi_R^k\to W>0$, there is $R_2\ge1$ with $\int w\phi_R^k\ge W/2$ for $R\ge R_2$. Hence, for $R\ge R_2$ such that $Z(R^2)\ge\norm{\zeta_-}_1$,
\begin{equation}\label{eq:nonint}
\frac W2\big(Z(R^2)-\norm{\zeta_-}_1\big)\le\norm{u_0}_1+2C_0R^{N-\frac2{p-1}} .
\end{equation}
In case (i), $N-2/(p-1)\le0$, so the right-hand side of \eqref{eq:nonint} is bounded while the left-hand side tends to $+\infty$. In case (ii), set $\beta:=\frac1{p-1}-\frac N2<0$ and choose $T_j=R_j^2\to\infty$ with $T_j^{\beta}Z(T_j)\to\infty$; in particular $Z(T_j)\to\infty$. Multiplying \eqref{eq:nonint} by $R_j^{2\beta}=R_j^{\frac2{p-1}-N}$ gives $\frac W2T_j^{\beta}Z(T_j)\le(\norm{u_0}_1+\frac W2\norm{\zeta_-}_1)T_j^{\beta}+2C_0$, whose right-hand side is bounded, a contradiction.
\end{proof}

\begin{remark}\rm \label{rem:power}
For $\zeta(t)=t^{\sigma}$, $\sigma>-1$, one has $Z(T)=T^{\sigma+1}/(\sigma+1)$, and Corollary~\ref{cor:nonintegrable} gives nonexistence whenever $\sigma+1>\frac N2-\frac1{p-1}$, that is for $p<\frac{N-2\sigma}{N-2\sigma-2}$ when $N>2\sigma+2$, and for every $p>1$ when $N\le2\sigma+2$. This recovers the subcritical blow-up range of Jleli, Kawakami and Samet \cite{JKS2020}, for weak solutions, in every dimension and for sign-changing $u_0\in L^1$; the critical case $p=p^*(\sigma)$ requires a finer argument and is treated in \cite{JKS2020}. As $\sigma\downarrow-1$ the exponent $\frac{N-2\sigma}{N-2\sigma-2}$ decreases to $\frac{N+2}{N}=\pF$. An integrable forcing such as $t^{\sigma}e^{-t}$ sits at this endpoint: its impulse is finite, the $T$-power in \eqref{eq:nonint} is $0$, and the blow-up range collapses exactly to $1<p\le\pF$, now subject to the effective-mass condition.

The logarithmic profile $\zeta(t)=(1+t)^{-1}$ lies exactly on the borderline between the two regimes. Here $Z(T)=\log(1+T)\to\infty$, so by Corollary~\ref{cor:nonintegrable}(i) there is no global solution for $1<p\le\pF$ as soon as $W>0$, \emph{irrespective of the size and sign of $u_0$}: the effective mass is $+\infty$. On the other hand, splitting the convolution at $t/2$ as in the proof of Lemma~\ref{lem:kernel} one finds $(\abs{\zeta}\conv\hN)(t)\le C\big(t^{-N/2}\log(1+t)+(1+t)^{-1}\int_0^{t/2}\hN\big)$ for $t\ge2$, which is $O(t^{-1/2}\log t)$ if $N=1$, $O(t^{-1}\log t)$ if $N=2$ and $O(t^{-1})$ if $N\ge3$. Hence $\omega_\zeta<\infty$ whenever $p>\max\{\pF,2\}$, and Corollary~\ref{cor:small} gives small-data global existence in that range. For $N\in\{1,2\}$, where $\pF\ge2$, the critical exponent for $\zeta=(1+t)^{-1}$ is therefore still $\pF$, but the dichotomy in the Fujita range no longer depends on the data. Note that Corollary~\ref{cor:small} does not require $\zeta\in L^1$.
\end{remark}

\section{Global existence}\label{sec:existence}

\subsection{The linear forced flow}
Let $u_0,w\in L^1(\RN)\cap L^\infty(\RN)$ and $\zeta\in L^1_{\rm loc}([0,\infty))$, and let $V$ and $m$ be defined by \eqref{eq:V-intro}. We record the elementary properties on which the proof of Theorem~\ref{thm:exist} rests. First, $V=S(t)\abs{u_0}+\Duh[\abs{\zeta}\otimes\abs{w}]$ is the solution of the linear problem with absolute data, and by Lemma~\ref{lem:heat}(a),(b),
\begin{equation}\label{eq:linear}
0\le V(x,t)\le m(t)\quad\text{for all }(x,t)\in\Q,
\qquad
S(t-s)V(s)\le V(t)\quad\text{for }0\le s\le t .
\end{equation}
The second inequality follows from the semigroup property, Tonelli's theorem and positivity:
$S(t-s)V(s)=S(t)\abs{u_0}+\int_0^s\abs{\zeta(r)}S(t-r)\abs{w}\,dr\le V(t)$. Second, by Lemma~\ref{lem:heat}(b),
\begin{equation}\label{eq:forcing}
m(t)\le\tn{u_0}\,\hN(t)+\tn{w}\,(\abs{\zeta}\conv\hN)(t)\qquad(t>0),
\end{equation}
which bounds the linear forced flow in terms of the forcing kernel. Third, the first inequality in \eqref{eq:linear} yields the pointwise interpolation inequality
\begin{equation}\label{eq:Lp-interp}
V(x,s)^{p}\le m(s)^{p-1}\,V(x,s),
\end{equation}
which replaces the usual $L^1$--$L^\infty$ interpolation $\norm{f}_p^p\le\norm{f}_1\norm{f}_\infty^{p-1}$ by a pointwise statement.

\subsection{Proof of Theorem~\ref{thm:exist}}
Let $0<T\le\infty$ and assume \eqref{eq:criterion}. Put
\[
a(t):=p\,\kappa^{p-1}m(t)^{p-1},\qquad A(t):=\int_0^ta(s)\ds ,
\]
and note that $A(t)\le p\kappa^{p-1}c_p=1$ on $(0,T)$, because $p\kappa^{p-1}c_p=p\cdot p^{p-1}(p-1)^{1-p}(p-1)^{p-1}p^{-p}=1$. Let $\mathcal K$ be the set of measurable functions $u$ on $\Q_T$ with $\abs{u}\le\kappa V$ a.e., where functions that agree a.e. are identified, and define
\[
d(u,v):=\operatorname*{ess\,sup}_{\{(x,t)\in\Q_T:\ V(x,t)>0\}}e^{-2A(t)}\,\frac{\abs{u(x,t)-v(x,t)}}{V(x,t)} .
\]
Every $u\in\mathcal K$ vanishes a.e. on $\{V=0\}$, so $d$ is a metric on $\mathcal K$ with $d\le2\kappa$. The space $(\mathcal K,d)$ is complete: if $(u_n)$ is $d$-Cauchy then, since $e^{-2A}\ge e^{-2}$, the quotients $u_n/V$ form a Cauchy sequence in $L^\infty(\{V>0\})$; their limit $\rho$ satisfies $\abs{\rho}\le\kappa$, and $u:=\rho V$ (extended by zero) belongs to $\mathcal K$ with $d(u_n,u)\to0$.

For $u\in\mathcal K$ we set
\[
\Phi[u](t):=S(t)u_0+\int_0^{t}S(t-s)\abs{u(s)}^{p}\ds+\int_0^tS(t-s)\big(\zeta(s)w\big)\ds .
\]
By \eqref{eq:Lp-interp}, $\abs{u}^p\le\kappa^pV^p\le\kappa^pm^{p-1}V$, and by \eqref{eq:linear},
$\int_0^tS(t-s)\big(m(s)^{p-1}V(s)\big)\ds\le\big(\int_0^tm^{p-1}\big)V(t)<\infty$, so $\Phi[u]$ is well defined and measurable.

\emph{$\Phi$ maps $\mathcal K$ into itself.} Using $\abs{S(\tau)f}\le S(\tau)\abs{f}$, the definition of $V$, \eqref{eq:Lp-interp}, \eqref{eq:linear} and \eqref{eq:criterion},
\begin{equation}\label{eq:nonlinear}
\begin{aligned}
\abs{\Phi[u](t)}&\le V(t)+\kappa^{p}\int_0^tS(t-s)V(s)^{p}\ds\\
&\le V(t)+\kappa^{p}\int_0^tm(s)^{p-1}S(t-s)V(s)\ds
\le\big(1+\kappa^{p}c_p\big)V(t)=\kappa\,V(t),
\end{aligned}
\end{equation}
since $\kappa^pc_p=p^p(p-1)^{-p}(p-1)^{p-1}p^{-p}=\frac1{p-1}$ and $1+\frac1{p-1}=\kappa$.

\emph{$\Phi$ is a contraction.} For $u,v\in\mathcal K$, the mean value theorem gives
\[
\big|\abs{u}^p-\abs{v}^p\big|\le p\max\{\abs u,\abs v\}^{p-1}\abs{u-v}\le a\,\abs{u-v},
\]
and by definition of $d$, $\abs{u-v}(s)\le d(u,v)\,e^{2A(s)}V(s)$ a.e. Hence, by \eqref{eq:linear},
\begin{equation}\label{eq:lip}
\begin{aligned}
\abs{\Phi[u]-\Phi[v]}(t)&\le\int_0^ta(s)\,S(t-s)\abs{u-v}(s)\ds\\
&\le d(u,v)\,V(t)\int_0^ta(s)e^{2A(s)}\ds
=d(u,v)\,V(t)\,\frac{e^{2A(t)}-1}{2},
\end{aligned}
\end{equation}
so that $d(\Phi[u],\Phi[v])\le\frac12d(u,v)$.

\emph{Conclusion.} By Banach's fixed point theorem $\Phi$ has a unique fixed point $u\in\mathcal K$. For every $T'<T$ we have $\abs{u}\le\kappa V\le\kappa\big(\norm{u_0}_\infty+\norm{w}_\infty\int_0^{T'}\abs{\zeta}\big)$ on $\Q_{T'}$ and $\norm{u(t)}_1\le\kappa\norm{V(t)}_1\le\kappa\big(\norm{u_0}_1+\norm{w}_1\int_0^{T'}\abs{\zeta}\big)$; hence $u$ is a mild solution on $(0,T)$ satisfying \eqref{eq:domination}. Uniqueness among all mild solutions follows from Lemma~\ref{lem:unique}. \qed

\begin{remark}\rm 
The constant $c_p$ is the best one allowed by the argument: a supersolution of the form $\kappa V$ requires $\kappa^p\int_0^Tm^{p-1}\le\kappa-1$, and $\max_{\kappa>1}(\kappa-1)\kappa^{-p}=(p-1)^{p-1}p^{-p}$ is attained at $\kappa=p/(p-1)$. If $u_0\not\equiv0$, then $m(t)\ge\norm{S(t)\abs{u_0}}_\infty\ge c\,t^{-N/2}$ for large $t$, and similarly if $\zeta\otimes w\not\equiv0$; therefore \eqref{eq:criterion} with $T=\infty$ forces $p>\pF$, in agreement with Theorem~\ref{thm:nonexist}.
\end{remark}

\subsection{Proof of Corollary~\ref{cor:small}}
Let $H_{N,p}:=\int_0^\infty\hN^{p-1}\dt$, which is finite since $p>\pF$ (Lemma~\ref{lem:heat}(f)). By \eqref{eq:forcing} and the elementary inequality $(\alpha+\beta)^{p-1}\le2^{(p-2)_+}(\alpha^{p-1}+\beta^{p-1})$ for $\alpha,\beta\ge0$,
\[
\int_0^\infty m^{p-1}\dt\le2^{(p-2)_+}\Big(H_{N,p}\tn{u_0}^{p-1}+\omega_\zeta\tn{w}^{p-1}\Big)
\le2^{(p-2)_+}\max\{1,H_{N,p}\}\Big(\tn{u_0}^{p-1}+\omega_\zeta\tn{w}^{p-1}\Big).
\]
With $\delta:=c_p\,2^{-(p-2)_+}/\max\{1,H_{N,p}\}$, the smallness assumption implies \eqref{eq:criterion} with $T=\infty$, and Theorem~\ref{thm:exist} applies. The bounds on $u(t)$ follow from $\abs{u}\le\kappa V$, \eqref{eq:linear} and \eqref{eq:forcing}, and from $\norm{V(t)}_1\le\norm{u_0}_1+\norm{w}_1\int_0^t\abs{\zeta}$.

It remains to verify $\omega_\zeta<\infty$ in cases (a)--(c). In case (a), Lemma~\ref{lem:kernel}(iii) with $q=p-1\ge1$ gives $\omega_\zeta^{1/(p-1)}\le\norm{\zeta}_1\norm{\hN}_{L^{p-1}(0,\infty)}<\infty$. In case (b), Lemma~\ref{lem:kernel}(i) gives (H$_K$), so $\omega_\zeta\le K^{p-1}H_{N,p}$, and \eqref{eq:forcing} yields $\norm{u(t)}_\infty\le\kappa(\tn{u_0}+K\tn{w})\hN(t)\le C(1+t)^{-N/2}$, because $\hN(t)\le2^{N/2}(1+4\pi t)^{-N/2}$. In case (c), Lemma~\ref{lem:kernel}(ii) gives $(\abs{\zeta}\conv\hN)^{p-1}\le C(1+t)^{-(p-1)\min\{b,N/2\}}$, which is integrable because $(p-1)b>1$ and $(p-1)N/2>1$. \qed

\begin{remark}\rm \label{rem:omega}
For $N\le2$ one has $\pF\ge2$, so case (a) covers every $p>\pF$ and every $\zeta\in L^1(0,\infty)$. For $N\ge3$ and $\pF<p<2$ the corollary requires some decay of $\zeta$, for instance $\abs{\zeta(t)}\lesssim t^{-b}$ with $b>1/(p-1)$. The reason is that the domination principle controls $\abs{u}^p$ through $\norm{V(t)}_\infty^{p-1}$, and the near-field part of the forced flow, of size $\abs{\zeta(t)}$ near the support of $w$, decays only as fast as $\zeta$ itself. Whether small-data global existence holds for every integrable $\zeta$ when $\pF<p<2$ is left open (Problem~\ref{prob:general-zeta}).
\end{remark}

\subsection{Proof of Theorem~\ref{thm:nonpositive}}
Set $v_0:=\abs{u_0}=-u_0$ and $g(x,t):=\abs{\zeta(t)w(x)}=-\zeta(t)w(x)\ge0$. Then $V=S(t)v_0+\Duh[g]$, and
\[
0\le V\le M_V:=\norm{u_0}_\infty+\norm{\zeta}_1\norm{w}_\infty\qquad\text{in }\Q .
\]
Let $\lambda:=pM_V^{p-1}$ and $\theta(r):=\lambda r-r^p$. On $[0,M_V]$ the function $\theta$ is nonnegative and nondecreasing, since $\theta(r)=r(\lambda-r^{p-1})$ and $\theta'(r)=\lambda-pr^{p-1}$. For measurable $v$ with $0\le v\le V$ we define
\[
\Phi[v](t):=e^{-\lambda t}S(t)v_0+\int_0^te^{-\lambda(t-s)}S(t-s)\big(g(s)+\theta(v(s))\big)\ds .
\]
Clearly $\Phi[v]\ge0$, and $\Phi$ is order preserving on the order interval $\{0\le v\le V\}$ because $\theta$ is nondecreasing on $[0,M_V]$ and $S(\tau)$ is positivity preserving. Moreover, since $\theta(V)\le\lambda V$, Lemma~\ref{lem:damped} applied on each $\Q_{T'}$ with $F=g$ (note that $\int_0^{T'}\norm{g(s)}_\infty\ds\le\norm{\zeta}_1\norm{w}_\infty$) gives
\[
\Phi[V](t)\le e^{-\lambda t}S(t)v_0+\int_0^te^{-\lambda(t-s)}S(t-s)\big(g(s)+\lambda V(s)\big)\ds=V(t).
\]
Thus $\Phi$ maps $\{0\le v\le V\}$ into itself. Define $v^0:=0$ and $v^{n+1}:=\Phi[v^n]$. By induction, $0\le v^n\le v^{n+1}\le V$ for all $n$: indeed $v^1=\Phi[0]\ge0=v^0$, and $v^{n}\le v^{n+1}$ implies $v^{n+1}=\Phi[v^n]\le\Phi[v^{n+1}]=v^{n+2}$. Let $v:=\lim_nv^n$. Since $\theta$ is continuous and nondecreasing on $[0,M_V]$, the integrands $g+\theta(v^n)$ are nonnegative and increase to $g+\theta(v)$, so the monotone convergence theorem gives $v=\Phi[v]$. As $v$ is bounded and $F:=g-v^p$ satisfies $\int_0^{T'}\norm{F(s)}_\infty\ds<\infty$, and since $g+\theta(v)=F+\lambda v$, the converse part of Lemma~\ref{lem:damped} yields
\[
v(t)=S(t)v_0+\int_0^tS(t-s)\big(g(s)-v(s)^p\big)\ds\qquad\text{a.e. in }\Q .
\]
Setting $u:=-v$ and recalling $-v_0=u_0$, $-g=\zeta w$ and $v^p=\abs{u}^p$, we obtain \eqref{eq:mild}. Since $-V\le u\le0$, $u$ is bounded on $\Q$ and $\norm{u(t)}_1\le\norm{V(t)}_1\le\norm{u_0}_1+\norm{\zeta}_1\norm{w}_1$; hence $u$ is a global mild solution, unique by Lemma~\ref{lem:unique}.

Assume finally $1<p\le\pF$. By Lemma~\ref{lem:mild-weak}(a), $u$ is a global weak solution, and Theorem~\ref{thm:annihilation} gives $\iQ\abs{u}^p=-\Meff=\norm{u_0}_1-ZW$. Since $\zeta w\le0$ we have $-ZW=-\iQ\zeta w=\norm{\zeta w}_{L^1(\Q)}$. The convergence $\norm{u(t)}_1\to0$ is a special case of Theorem~\ref{thm:asymptotics}(ii), proved in Section~\ref{sec:asymptotics} independently of the present theorem. \qed

\section{Sharp lifespan in the Fujita range}\label{sec:lifespan}

Throughout this section $1<p\le\pF$, $\zeta\in L^1(0,\infty)$, $a,b\in L^1(\RN)\cap L^\infty(\RN)$, and for $\eps>0$ we consider \eqref{eq:model} with $u_0=\eps a$ and $w=\eps b$. By Lemma~\ref{lem:unique}, mild solutions on different intervals coincide on their common interval, so $T_\eps$ is the maximal existence time of the mild solution.

\subsection{Upper bounds}
Assume $\Meff_*=\int a+Z\int b>0$, let $\eps\in(0,1]$, and suppose that \eqref{eq:model} has a mild solution $u$ on $(0,T)$. By Lemma~\ref{lem:mild-weak}(a), $u$ is a weak solution on $\Q_T$. For $0<R<R_T:=\sqrt{T/2}$ the parabolic cutoff $\Psi_R$ is supported in $\{\abs{x}^2+t\le2R^2\}$, a compact subset of $\RN\times[0,T)$, and is therefore an admissible test function on $\Q_T$. Exactly as in \eqref{eq:identity} we obtain
\begin{equation}\label{eq:identity-eps}
J(R)+\eps\,c(R)=-\iint_{\Q_T}u\,\Lop\Psi_R,\qquad
c(R):=\int_{\RN}a\,\eta\big(\abs{x}^2/R^2\big)^k\dx+\iQ\zeta\,b\,\Psi_R ,
\end{equation}
where $J(R)=\iint_{\Q_T}\abs{u}^p\Psi_R$. By dominated convergence $c(R)\to\Meff_*$ as $R\to\infty$, so there is $R_1\ge1$, depending only on $a$, $b$, $\zeta$ and $\eta$, such that $c(R)\ge\Meff_*/2$ for all $R\ge R_1$. If $T\le2R_1^2$ there is nothing to prove, so we assume $R_T>R_1$.

\emph{The subcritical case $p<\pF$.} The argument leading to \eqref{eq:subcritical} applies verbatim on $\Q_T$ and gives, for $R_1\le R<R_T$,
\[
\frac{\eps\Meff_*}2\le\eps\,c(R)\le\frac{C_2^{p'}}{p'}\,R^{-(2p'-N-2)} .
\]
Letting $R\uparrow R_T$ we get $T\le2\big(2C_2^{p'}/(p'\Meff_*\eps)\big)^{2/(2p'-N-2)}$, and $\frac{2}{2p'-N-2}=\frac{2(p-1)}{2-N(p-1)}$ because $2p'-2=\frac2{p-1}$. Hence $T\le\max\{2R_1^2,C\eps^{-2(p-1)/(2-N(p-1))}\}\le C'\eps^{-2(p-1)/(2-N(p-1))}$ for $\eps\in(0,1]$, and the same bound holds for $T_\eps$.

\emph{The critical case $p=\pF$.} Now $N+2-2p'=0$ and \eqref{eq:Holder-par} reads
\begin{equation}\label{eq:crit-ineq}
J(R)+\eps\,c(R)\le C_2\,y(R)^{1/p},\qquad y(\rho):=\iint_{\Q_T}\abs{u}^p\,\eta_*\big(s_\rho\big)\dx\dt,\qquad \eta_*:=\eta^{k}\mathbf 1_{[1,2]},
\end{equation}
because $\Psi_R\mathbf 1_{\mathcal D_R}=\eta_*(s_R)$. Following the idea of Ikeda and Sobajima \cite{IkedaSobajima2019}, we introduce
\[
Y(R):=\int_0^R y(\rho)\,\frac{d\rho}{\rho},\qquad 0<R<R_T .
\]
By Tonelli's theorem, $Y(R)=\iint_{\Q_T}\abs{u}^p\,\Theta_R(x,t)\dx\dt$ with $\Theta_R(x,t):=\int_0^R\eta_*\big((\abs{x}^2+t)/\rho^2\big)\frac{d\rho}{\rho}$. For $r:=\abs{x}^2+t>0$ the substitution $\varsigma=r/\rho^2$, $\frac{d\rho}{\rho}=-\frac{d\varsigma}{2\varsigma}$, gives
\[
\Theta_R(x,t)=\frac12\int_{r/R^2}^{\infty}\eta_*(\varsigma)\frac{d\varsigma}{\varsigma}
\le\frac12\,\eta\big(r/R^2\big)^{k}\int_1^2\frac{d\varsigma}{\varsigma}=\frac{\log2}{2}\,\Psi_R(x,t),
\]
where we used that $\eta_*$ vanishes outside $[1,2]$ and that $\eta(\varsigma)\le\eta(r/R^2)$ for $\varsigma\ge r/R^2$. Consequently
\begin{equation}\label{eq:Y-J}
Y(R)\le\frac{\log2}{2}\,J(R)<\infty .
\end{equation}
In particular $y(\rho)/\rho$ is integrable on $(0,R]$ for every $R<R_T$, so $Y$ is absolutely continuous on $[R_1,R]$ with $Y'(R)=y(R)/R$ for a.e. $R$. Combining \eqref{eq:crit-ineq}, \eqref{eq:Y-J} and $c(R)\ge\Meff_*/2$, we find for a.e. $R\in[R_1,R_T)$
\[
\frac{2}{\log2}\,\widetilde Y(R)\le C_2\big(R\,\widetilde Y'(R)\big)^{1/p},\qquad
\widetilde Y(R):=Y(R)+\frac{\log2}{4}\,\Meff_*\,\eps>0 .
\]
Thus $\widetilde Y'(R)\ge c_3\widetilde Y(R)^p/R$ with $c_3:=\big(2/(C_2\log2)\big)^p$. Since $\widetilde Y$ is absolutely continuous and bounded below by a positive constant, $\widetilde Y^{1-p}$ is absolutely continuous and $(\widetilde Y^{1-p})'=(1-p)\widetilde Y^{-p}\widetilde Y'\le-(p-1)c_3/R$ a.e. Integrating over $[R_1,R]$ and using $\widetilde Y(R)^{1-p}>0$,
\[
(p-1)\,c_3\log\frac{R}{R_1}<\widetilde Y(R_1)^{1-p}\le\Big(\frac{\log2}{4}\Meff_*\eps\Big)^{1-p} .
\]
Letting $R\uparrow R_T$ we obtain $T\le2R_1^2\exp\big(2C_4\eps^{-(p-1)}\big)$ with $C_4:=\big(\frac{\log2}{4}\Meff_*\big)^{1-p}/((p-1)c_3)$, and therefore $T_\eps\le\exp(C\eps^{-(p-1)})$ for $\eps\in(0,1]$, with $C:=2C_4+\log_+(2R_1^2)$. This proves Theorem~\ref{thm:lifespan}(i). \qed

\subsection{Lower bounds}
Assume (H$_K$) (by Lemma~\ref{lem:kernel}(i) this is implied by (H$_a$)), and set $\Lambda:=\tn{a}+K\tn{b}$. If $\Lambda=0$ the data vanish and $T_\eps=\infty$, so let $\Lambda>0$. By \eqref{eq:forcing} with $u_0=\eps a$, $w=\eps b$,
\[
m(t)\le\eps\big(\tn{a}\hN(t)+\tn{b}(\abs{\zeta}\conv\hN)(t)\big)\le\eps\Lambda\,\hN(t),
\]
hence
\[
\int_0^Tm^{p-1}\dt\le(\eps\Lambda)^{p-1}H(T),\quad H(T):=\int_0^T\hN^{p-1}\dt .
\]
By Theorem~\ref{thm:exist}, $T_\eps\ge T$ whenever $(\eps\Lambda)^{p-1}H(T)\le c_p$. Let $\beta:=N(p-1)/2\le1$.

If $p<\pF$, then $\beta<1$ and $H(T)\le\int_0^T(4\pi t)^{-\beta}\dt=\frac{(4\pi)^{-\beta}}{1-\beta}T^{1-\beta}$. Choosing $T=\big(c_p(1-\beta)(4\pi)^{\beta}\big)^{1/(1-\beta)}(\eps\Lambda)^{-(p-1)/(1-\beta)}$ we obtain $T_\eps\ge c\,\eps^{-(p-1)/(1-\beta)}$, and $\frac{p-1}{1-\beta}=\frac{2(p-1)}{2-N(p-1)}$.

If $p=\pF$, then $\beta=1$ and $H(T)=\frac1{4\pi}\big(1+\log(4\pi T)\big)$ for $T\ge\frac1{4\pi}$. Choosing
\[
T=\tfrac1{4\pi}\exp\big(4\pi c_p(\eps\Lambda)^{-(p-1)}-1\big),
\]
which is admissible as soon as $4\pi c_p(\eps\Lambda)^{-(p-1)}\ge1$, we obtain
\[
T_\eps\ge\tfrac1{4\pi}\exp\big(4\pi c_p\Lambda^{1-p}\eps^{-(p-1)}-1\big)\ge\exp\big(2\pi c_p\Lambda^{1-p}\eps^{-(p-1)}\big)
\]
for $\eps\le\eps_0$ small enough. This proves Theorem~\ref{thm:lifespan}(ii). \qed

\begin{remark}\rm \label{rem:lifespan}
The upper bounds use only $\zeta\in L^1(0,\infty)$ and $a,b\in L^1$, and hold for the lifespan of weak solutions; the lower bounds hold for arbitrary signs. Taking $a=0$ shows that a small forcing alone, with $Z\int b>0$, produces blow-up with the same lifespan scale as a small initial datum of positive mass. The constants depend on the data through $\Meff_*$ (upper bounds) and $\Lambda$ (lower bounds); for nonnegative data and $\zeta\ge0$ they are comparable. The lifespan of small data with $\Meff_*\le0$ is not addressed here (Problem~\ref{prob:lifespan}).
\end{remark}

\section{Large-time behavior and asymptotic mass}\label{sec:asymptotics}

\begin{lemma}\label{lem:duhamel-asymp}
Let $\varphi\in L^1(\RN)$ and $f\in L^1(\Q)$. Then
\[
\lim_{t\to\infty}\Big\|S(t)\varphi+\int_0^tS(t-s)f(s)\ds-\Big(\int_{\RN}\varphi\dx+\iQ f\dx\ds\Big)G_t\Big\|_1=0 .
\]
\end{lemma}

\begin{proof}
The contribution of $\varphi$ is handled by Lemma~\ref{lem:heat}(d). For the Duhamel term we write
\[
\int_0^tS(t-s)f(s)\ds-\Big(\iQ f\Big)G_t
=\int_0^t\Big[S(t-s)f(s)-\Big(\int f(s)\Big)G_t\Big]\ds-\Big(\int_t^\infty\!\!\int_{\RN}f\dx\ds\Big)G_t .
\]
The $L^1$ norm of the last term is at most $\int_t^\infty\norm{f(s)}_1\ds\to0$. For the first term, set $\delta(t,s):=\mathbf 1_{\{s<t\}}\norm{S(t-s)f(s)-(\int f(s))G_t}_1$. Then $\delta(t,s)\le2\norm{f(s)}_1$, which is integrable in $s$, and for a.e. fixed $s$,
\[
\delta(t,s)\le\Big\|S(t-s)f(s)-\Big(\int f(s)\Big)G_{t-s}\Big\|_1+\norm{f(s)}_1\norm{G_{t-s}-G_t}_1\xrightarrow[t\to\infty]{}0
\]
by Lemma~\ref{lem:heat}(d),(e). Dominated convergence gives $\int_0^\infty\delta(t,s)\ds\to0$.
\end{proof}

\subsection{Proof of Theorem~\ref{thm:asymptotics}}
(i) Let $u$ be a global mild solution with $\abs{u}^p\in L^1(\Q)$. Then $f:=\abs{u}^p+\zeta w\in L^1(\Q)$, since $\norm{\zeta\otimes w}_{L^1(\Q)}=\norm{\zeta}_1\norm{w}_1$, and $\iQ f=\iQ\abs{u}^p+ZW$. Applying Lemma~\ref{lem:duhamel-asymp} to \eqref{eq:mild} with $\varphi=u_0$ gives $\norm{u(t)-m_\infty G_t}_1\to0$ with $m_\infty=m_0+ZW+\iQ\abs{u}^p=\Meff+\iQ\abs{u}^p$.

(ii) Let $1<p\le\pF$ and let $u$ be a global mild solution. By Lemma~\ref{lem:mild-weak}(a) it is a global weak solution, and Theorem~\ref{thm:annihilation} gives $\iQ\abs{u}^p=-\Meff<\infty$. Part (i) applies with $m_\infty=\Meff+\iQ\abs{u}^p=0$, so $\norm{u(t)}_1=\norm{u(t)-0\cdot G_t}_1\to0$.

(iii) Let $p>\pF$, let $\zeta$ satisfy (H$_a$), and let $u$ be a global mild solution with $\abs{u}\le\kappa V$. By Lemma~\ref{lem:kernel}(i), (H$_K$) holds, so that \eqref{eq:linear} and \eqref{eq:forcing} give
\[
\norm{u(s)}_\infty\le\kappa\,m(s)\le\kappa\big(\tn{u_0}+K\tn{w}\big)\hN(s)\le C_5(1+s)^{-N/2},\]
\[
\norm{u(s)}_1\le\kappa\big(\norm{u_0}_1+\norm{\zeta}_1\norm{w}_1\big)=:C_6 .
\]
Hence $\iQ\abs{u}^p\le\int_0^\infty\norm{u(s)}_\infty^{p-1}\norm{u(s)}_1\ds\le C_5^{p-1}C_6\int_0^\infty(1+s)^{-N(p-1)/2}\ds<\infty$, because $N(p-1)/2>1$. Part (i) gives the $L^1$ convergence. For the $L^\infty$ estimate, let $t\ge2$. By Lemma~\ref{lem:mild-weak}(c) with $\tau=t/2$ and the identity $S(t/2)G_{t/2}=G_t$,
\[
\begin{aligned}
u(t)-m_\infty G_t={}&S(t/2)\big[u(t/2)-m_\infty G_{t/2}\big]\\&+\int_{t/2}^tS(t-s)\abs{u(s)}^p\ds+\int_{t/2}^tS(t-s)\big(\zeta(s)w\big)\ds .
\end{aligned}
\]
By Lemma~\ref{lem:heat}(b) and part (i), the first term has $L^\infty$ norm at most $(2\pi t)^{-N/2}\norm{u(t/2)-m_\infty G_{t/2}}_1=o(t^{-N/2})$. The second term is bounded in $L^\infty$ by
\[
\int_{t/2}^t\norm{u(s)}_\infty^p\ds\le\frac t2\,C_5^p\,(1+t/2)^{-Np/2}=O(t^{1-Np/2})=o(t^{-N/2}),
\]
again because $N(p-1)/2>1$. The third term is bounded by $\tn{w}\int_{t/2}^t\abs{\zeta(s)}\hN(t-s)\ds=o(t^{-N/2})$ by \eqref{eq:o-cond}. Hence $t^{N/2}\norm{u(t)-m_\infty G_t}_\infty\to0$, and the case $1<q<\infty$ follows from the inequality $\norm{g}_q\le\norm{g}_1^{1/q}\norm{g}_\infty^{1-1/q}$:
\[
t^{\frac N2(1-\frac1q)}\norm{u(t)-m_\infty G_t}_q\le\norm{u(t)-m_\infty G_t}_1^{1/q}\Big(t^{N/2}\norm{u(t)-m_\infty G_t}_\infty\Big)^{1-1/q}\xrightarrow[t\to\infty]{}0 .
\]
Finally $m_\infty-\Meff=\iQ\abs{u}^p\ge0$. \qed

\begin{corollary}\label{cor:characterization}
Let $\zeta$ satisfy \textup{(H$_a$)} and let $w\in L^1(\RN)\cap L^\infty(\RN)$. Consider \eqref{eq:model} with forcing $\mu\,\zeta(t)w(x)$, $\mu\in\R$.
\begin{enumerate}
\item[(i)] If $1<p\le\pF$, then for every $\mu\in\R$ and every $u_0\in L^1\cap L^\infty$, every global mild solution satisfies $\norm{u(t)}_1\to0$ as $t\to\infty$.
\item[(ii)] If $p>\pF$, there exist $\mu\ne0$ and $u_0\in L^1\cap L^\infty$ for which \eqref{eq:model} has a global mild solution with $\liminf_{t\to\infty}\norm{u(t)}_1>0$.
\end{enumerate}
Thus $p\le\pF$ if and only if all global solutions lose their mass.
\end{corollary}

\begin{proof}
Part (i) is Theorem~\ref{thm:asymptotics}(ii). For (ii), fix $a\in L^1\cap L^\infty$ with $a\ge0$ and $\int a>\abs{Z W}$, and take $u_0=\eps a$, $\mu=\eps$. Then $\Meff=\eps(\int a+ZW)>0$. For $\eps$ small, Corollary~\ref{cor:small}(b) provides a global mild solution with $\abs{u}\le\kappa V$, and by Theorem~\ref{thm:asymptotics}(iii) and Lemma~\ref{lem:heat}(c), $\norm{u(t)}_1\ge\abs{\int u(t)}\to m_\infty\ge\Meff>0$.
\end{proof}

\begin{remark}\rm 
For $1<p\le\pF$, the vanishing of $\norm{u(t)}_1$ holds for sign-changing solutions and in the presence of forcing, and it does not rely on any comparison argument: it is a consequence of the identity \eqref{eq:annihilation} and of the linear asymptotics of Lemma~\ref{lem:duhamel-asymp}. For $p>\pF$ the asymptotic mass $m_\infty=\Meff+\iQ\abs{u}^p$ exceeds the effective mass by the total production of the source. For the small solutions of Corollary~\ref{cor:small}, $\iQ\abs{u}^p\le\kappa^p\big(\int_0^\infty m^{p-1}\big)\sup_t\norm{V(t)}_1$ by \eqref{eq:Lp-interp}, so that $m_\infty=\Meff+O(\eps^{p})$ when the data are of size $\eps$.
\end{remark}

\section{The model example \texorpdfstring{$\zeta(t)=t^{\sigma}e^{-t}$}{zeta(t)=t^sigma e^{-t}}}\label{sec:model}

We now collect the results for the model forcing \eqref{eq:model-zeta}, for which $Z=\Gamma(\sigma+1)$ and (H$_a$) holds with every $a>0$.

\begin{theorem}\label{thm:model}
Let $\zeta(t)=t^{\sigma}e^{-t}$ with $\sigma>-1$, let $u_0,w\in L^1(\RN)\cap L^\infty(\RN)$, and let $\Meff=m_0+\Gamma(\sigma+1)W$.
\begin{enumerate}
\item[(i)] \textup{(Fujita range.)} If $1<p\le\pF$ and \eqref{eq:model} has a global weak solution, then either $\Meff<0$, or $u_0\equiv0$ and $w\equiv0$. Every global weak solution satisfies $\iQ\abs{u}^p=-\Meff$, and every global mild solution satisfies $\norm{u(t)}_1\to0$. When $\Meff<0$, global solutions and blowing-up solutions both occur.
\item[(ii)] \textup{(Above the Fujita exponent.)} If $p>\pF$, there exist $\eps>0$ and $C_*>0$, depending only on $N$, $p$ and $\sigma$, such that if $\tn{u_0}+\tn{w}\le\eps$, then \eqref{eq:model} possesses a unique global mild solution. It satisfies, for all $t>0$,
\[
\norm{u(t)}_{L^1}\le C_\ast\eps,
\qquad
\norm{u(t)}_{L^{\infty}}\le C_\ast\eps\,(1+t)^{-N/2},
\]
and, with $m_\infty:=\Meff+\iQ\abs{u}^p$, one has $0\le m_\infty-\Meff\le C_*\eps^{p}$ and, for every $q\in[1,\infty]$, $t^{\frac N2(1-\frac1q)}\norm{u(t)-m_\infty G_t}_q\to0$. No sign restriction on $u_0$ or $w$ is required.
\item[(iii)] \textup{(Lifespan.)} If $1<p\le\pF$, $u_0=\eps a$, $w=\eps b$ and $\int a+\Gamma(\sigma+1)\int b>0$, then for $\eps$ small $T_\eps\simeq\eps^{-2(p-1)/(2-N(p-1))}$ if $p<\pF$, and $\log T_\eps\simeq\eps^{-(p-1)}$ if $p=\pF$.
\item[(iv)] \textup{(Large data.)} For every $p>1$, data $u_0=\lambda g$ with $\abs{\{g>0\}}>0$ give no global weak solution when $\lambda$ is large, whatever $w$.
\end{enumerate}
\end{theorem}

\begin{proof}
Since $\zeta\in L^1(0,\infty)$ and $\abs{\zeta(t)}\le C_at^{-a}$ for $t\ge1$ and every $a>0$, (H$_a$) holds, and by Lemma~\ref{lem:kernel}(i) so does (H$_K$) with a constant $K=K(N,\sigma)$. Part (i) follows from Theorems~\ref{thm:nonexist} and \ref{thm:asymptotics}(ii) and Example~\ref{ex:both}. In (ii), $\omega_\zeta\le K^{p-1}H_{N,p}$ depends only on $N,p,\sigma$, and $\tn{u_0}^{p-1}+\omega_\zeta\tn{w}^{p-1}\le(1+\omega_\zeta)\eps^{p-1}\le\delta$ for $\eps$ small; Corollary~\ref{cor:small}(b) gives the solution and the bounds, with $\norm{u(t)}_1\le\kappa(\norm{u_0}_1+\Gamma(\sigma+1)\norm{w}_1)$. The bound on $m_\infty-\Meff$ follows from the last remark of Section~\ref{sec:asymptotics}, since $\int_0^\infty m^{p-1}\le C\eps^{p-1}$ and $\sup_t\norm{V(t)}_1\le C\eps$, and the asymptotics from Theorem~\ref{thm:asymptotics}(iii). Parts (iii) and (iv) are Theorem~\ref{thm:lifespan} and Corollary~\ref{cor:obstruction}(a).
\end{proof}

\begin{remark}\rm \label{rem:correction}
For growing forcings $t^{\sigma}w$, $\sigma>0$, there is no global solution for any $p>1$ when $W>0$ \cite{JKS2020}, and one might expect a similar conclusion for profiles such as $t^{\sigma}e^{-t}$ which are large on bounded time intervals. This is not the case. The exponential factor makes $\zeta\in L^1(0,\infty)$, so the forcing delivers the \emph{finite} total impulse $\Gamma(\sigma+1)W$ and the lower bound \eqref{eq:F+I-limit} saturates at the constant $\Meff$ instead of growing with $T$. In the test-function balance this $O(1)$ gain can defeat the diffusion cost $R^{\,N+2-2p'}$ only when $N+2-2p'\le0$, i.e.\ for $p\le\pF$. Theorem~\ref{thm:model}(ii) shows that the restriction is sharp: for $p>\pF$ small data yield global solutions, so nonexistence \emph{cannot} hold for every $p>1$. The same conclusion applies to every forcing profile $\zeta$ that is integrable and decays faster than any power, since such a $\zeta$ satisfies (H$_a$) for every $a$; the blow-up range for such forcings is therefore $1<p\le\pF$ (under the effective-mass condition), not $p>1$. What remains true for every $p>1$ is blow-up for \emph{large} data or forcing (Corollary~\ref{cor:obstruction}).
\end{remark}

\section{Concluding remarks and open problems}\label{sec:open}

The results of Sections~\ref{sec:annihilation}--\ref{sec:model} are proved in full. In this section we collect questions that we were not able to settle; everything stated here is conjectural or open and is clearly separated from the results above.

\begin{remark}\rm \label{rem:open}
Theorem~\ref{thm:exist} and Corollary~\ref{cor:small} are small-data statements. For $p>\pF$ large data blow up (Corollary~\ref{cor:obstruction}(a),(b)); because the nonlinearity is $\abs{u}^{p}$, the source is nonnegative and only the positive part of the data drives blow-up, while negative data are damped (Theorem~\ref{thm:nonpositive}). In the preliminary version of this work it was asked whether the hypothesis \eqref{eq:mass-condition} in Theorem~\ref{thm:nonexist} could be removed altogether for $1<p\le\pF$, that is, whether a very negative effective mass still forces blow-up through the $\abs{u}^{p}$ source. Theorem~\ref{thm:nonpositive} answers this in the negative, and Theorem~\ref{thm:nonexist}(ii) shows that the condition can be relaxed exactly to $\Meff\ge0$.
\end{remark}

\begin{problem}\label{prob:negative}
Let $1<p\le\pF$. Describe the set of data $(u_0,\zeta w)$ with $\Meff<0$ for which \eqref{eq:model} has a global solution. In particular, is the nonpositive global solution of Theorem~\ref{thm:nonpositive} stable under small perturbations of the data that keep $\Meff<0$, for instance in the norm $\tn{\cdot}$? By \eqref{eq:annihilation} any such global solution produces exactly $\abs{\Meff}$ units of mass through the source, which suggests that the answer depends on how the positive part of the perturbation is distributed relative to the negative mass.
\end{problem}

\begin{problem}\label{prob:lifespan}
Let $1<p\le\pF$ and $u_0=\eps a$, $w=\eps b$ with $\Meff_*=\int a+Z\int b\le0$ and $(a,\zeta b)$ not of the sign considered in Theorem~\ref{thm:nonpositive}. Theorem~\ref{thm:lifespan}(ii) gives the lower bounds in all cases. When $\Meff_*=0$ and $a$ has a nonzero first moment, one expects the lifespan to be governed by a higher-order term of the asymptotic expansion of the linear flow, as for the unforced equation with zero-mass data. Determine the sharp lifespan in this case.
\end{problem}

\begin{problem}\label{prob:general-zeta}
Let $N\ge3$ and $\pF<p<2$. Does small-data global existence hold for \emph{every} $\zeta\in L^1(0,\infty)$? Remark~\ref{rem:omega} explains why the domination principle requires some decay of $\zeta$. A heuristic computation of $\iQ F^p$, where $F$ is the linear forced flow, suggests that the answer is positive, but a proof would require a finer description of the near-field part of $F$ than the $L^\infty$ bound \eqref{eq:forcing}.
\end{problem}

\begin{problem}\label{prob:rv}
For $\zeta$ regularly varying at infinity with index $-1$, so that $Z(T)=\int_0^T\zeta$ is slowly varying in the sense of Karamata, Corollary~\ref{cor:nonintegrable} and Remark~\ref{rem:power} show that the critical exponent is $\pF$ (at least for $N\le2$) and that the Fujita range no longer depends on the data when $Z(T)\to\infty$ and $W>0$. For profiles of index $\sigma\in(-1,0)$ with slowly varying corrections, determine the behavior at the threshold $p=p^*(\sigma)$, where $T^{\frac1{p-1}-\frac N2}Z(T)$ is itself slowly varying; one expects blow-up if and only if this quantity is unbounded, in the spirit of \cite{KumarMajdoub2026}.
\end{problem}

\begin{problem}\label{prob:rates}
For $1<p\le\pF$ every global mild solution satisfies $\norm{u(t)}_1\to0$. For the absorption equation with nonnegative data the rate is known: $\norm{v(t)}_1\simeq t^{-(\frac1{p-1}-\frac N2)}$ for $p<\pF$ \cite{BrezisPeletierTerman1986,KaminPeletier1985} and $\norm{v(t)}_1\simeq(\log t)^{-N/2}$ for $p=\pF$ \cite{GalaktionovKurdyumovSamarskii1986}. Do sign-changing global solutions of \eqref{eq:model} with integrable forcing decay at the same rates, possibly under a decay assumption on $\zeta$?
\end{problem}

\begin{problem}\label{prob:other}
The proof of Theorem~\ref{thm:annihilation} uses the operator only through the scaling of the cutoff cost \eqref{eq:cost-par}. We expect the identity \eqref{eq:annihilation}, with $\pF$ replaced by the corresponding critical exponent, to hold for the fractional Laplacian (with $\pF=1+2s/N$), for mixed local--nonlocal operators \cite{KiraneFinoAyoub2025,BenBelgacemMajdoub2025b}, for time-fractional diffusions \cite{BenBelgacemMajdoub2025c,FinoKirane2012}, and for diffusions with interfaces \cite{MajdoubMliki2026}. For nonlocal operators the cutoffs must have polynomial tails, and the proof of the vanishing of the diffusive term \eqref{eq:tail} has to be adapted accordingly.
\end{problem}

\bigskip

\section*{Declarations}
\noindent\textbf{Funding.} This research received no external funding.

\noindent\textbf{Conflict of interest.} The authors declare that they have no conflicts of interest.

\noindent\textbf{Data availability.} No data were generated or analyzed in this study.


\end{document}